\documentclass[a4paper,12pt]{amsart}
\usepackage{amssymb}
\RequirePackage{calrsfs}

\DeclareSymbolFont{rsfscript}{OMS}{rsfs}{m}{b}
\DeclareSymbolFontAlphabet{\mathrsfs}{rsfscript}
\renewcommand{\mathcal}{\mathrsfs}

\def\bfit{\bfseries\itshape}

\input xypic
\xyoption{all}
\xyoption{arc}

\newtheorem{theo}{Th\'eor\`eme}[section]

\newtheorem{prop}[theo]{Proposition}

\newtheorem{lem}[theo]{Lemme}

\newtheorem{coro}[theo]{Corollaire}

\def\remark#1{{\refstepcounter{theo}\label{#1}\noindent\sc Remarque
\arabic{section}.\arabic{theo} - }}
\def\example#1{{\refstepcounter{theo}\label{#1}\noindent\sc Exemple 
\arabic{section}.\arabic{theo} - }}

\def\equat{\refstepcounter{theo}$$~}
\def\endequat{\leqno{\boldsymbol{(\arabic{section}.\arabic{theo})}}~$$}

\def\AG{{\mathfrak A}}

    \def\NM{{\mathbb{N}}}

    \def\QM{{\mathbb{Q}}}
    \def\RM{{\mathbb{R}}}
\def\SG{{\mathfrak S}}

    \def\ZM{{\mathbb{Z}}}

\def\Bb{{\mathbf B}}    \def\BC{{\mathcal{B}}}
    \def\CC{{\mathcal{C}}}
    \def\DC{{\mathcal{D}}}
    \def\EC{{\mathcal{E}}}
    \def\FC{{\mathcal{F}}}
    
    \def\HC{{\mathcal{H}}}

    \def\LC{{\mathcal{L}}}

    \def\OC{{\mathcal{O}}}
    \def\PC{{\mathcal{P}}}
    
    \def\RC{{\mathcal{R}}}
    \def\SC{{\mathcal{S}}}
    
    \def\UC{{\mathcal{U}}}
    \def\VC{{\mathcal{V}}}

\def\Brm{{\mathrm{B}}}

    \def\HCB{{\boldsymbol{\mathcal{H}}}}

\def\Fti{{\tilde{F}}}

 \def\HCt{{\widetilde{\mathcal{H}}}}
    
\def\Jti{{\widetilde{J}}}

  \def\tti{{\tilde{t}}}

\def\Wti{{\widetilde{W}}}

\def\Sba{{\bar{S}}}          \def\sba{{\bar{s}}}

\def\Cov{{\overline{C}}}

          \def\hov{{\overline{h}}}

\def\g{\gamma}
\def\G{\Gamma}

\def\e{\varepsilon}
\def\ph{\varphi}
\def\l{\lambda}
\def\L{\Lambda}
\def\o{\omega}

\def\r{\rho}
\def\s{\sigma}

\def\th{\theta}

\def\t{\tau}

\def\phb{{\boldsymbol{\varphi}}}        \def\pht{{\tilde{\varphi}}}

          \def\psit{{\tilde{\psi}}}

\def\phba{{\bar{\varphi}}}

\def\psiba{{\bar{\psi}}}

\def\sigba{{\bar{\s}}}

\DeclareMathOperator{\Hom}{{\mathrm{Hom}}}

\DeclareMathOperator{\Id}{{\mathrm{Id}}}

\DeclareMathOperator{\Ker}{{\mathrm{Ker}}}

\DeclareMathOperator{\Pos}{{\mathrm{Pos}}}
\DeclareMathOperator{\Posbar}{{\overline{\mathrm{Pos}}}}

\def\imp{\Rightarrow}

\def\to{\rightarrow}
\def\longto{\longrightarrow}

\def\fonction#1#2#3#4#5{\begin{array}{rccc}
{#1} : & {#2} & \longto & {#3} \\
& {#4} & \longmapsto & {#5} 
\end{array}}

\def\fonctio#1#2#3#4{\begin{array}{ccc}
{#1} & \longto & {#2} \\
{#3} & \longmapsto & {#4} 
\end{array}}

\def\vide{\varnothing}

\def\DS{\displaystyle}
\def\SS{\scriptstyle}

\def\finl{~$\SS \square$}

\def\infspe{\hspace{0.1em}\mathop{\preccurlyeq}\nolimits\hspace{0.1em}}

\def\lexp#1#2{\kern\scriptspace\vphantom{#2}^{#1}\kern-\scriptspace#2}
\def\le{\mathop{\leqslant}\nolimits}
\def\ge{\mathop{\geqslant}\nolimits}

\mathchardef\inferieur="321E
\mathchardef\superieur="321F

\def\eqna{\begin{eqnarray*}}
\def\endeqna{\end{eqnarray*}}

\def\tors{tores maximaux }

\def\tors{{\mathrm{tors}}}

\def\itemth#1{\item[${\mathrm{(#1)}}$]}

\catcode`\@=11
\long\def\@car#1#2\@nil{#1}
\long\def\@first#1#2{#1}
\long\def\@second#1#2{#2}
\long\def\ifempty#1{\expandafter\ifx\@car#1@\@nil @\@empty
  \expandafter\@first\else\expandafter\@second\fi}
\catcode`\@=12

\def\positif{\PC\! os}

\DeclareMathOperator{\can}{{\mathrm{can}}}
\DeclareMathOperator{\Cell}{{\mathrm{Cell}}}
\DeclareMathOperator{\maps}{{\mathrm{\FC}}}
\DeclareMathOperator{\weight}{{\mathrm{Poids}}}
\DeclareMathOperator{\codim}{{\mathrm{codim}}}

\def\ellb{\boldsymbol{\ell}}
\def\ellt{\tilde{\ell}}

\def\aug{{\mathrm{aug}}}
\def\espace{\vphantom{\DS{\frac{A}{A}}}}
\def\dotcup{\hskip1mm\dot{\cup}\hskip1mm}

\begin{document}

\baselineskip=16pt

\title{Semi-continuit\'e des cellules de Kazhdan-Lusztig}

\author{C\'edric Bonnaf\'e}
\address{\noindent 
Labo. de Math. de Besan\c{c}on (CNRS: UMR 6623), 
Universit\'e de Franche-Comt\'e, 16 Route de Gray, 25030 Besan\c{c}on
Cedex, France} 

\makeatletter
\email{cedric.bonnafe@univ-fcomte.fr}

\subjclass{According to the 2000 classification:
Primary 20C08; Secondary 20C15}

\date{\today}
\def\abstractname{R\'esum\'e}

\begin{abstract} 
Des calculs dans les petits groupes de Coxeter et les groupes 
di\'edraux infinis 
sugg\`erent que les cellules de Kazhdan-Lusztig \`a param\`etres 
in\'egaux ob\'eissent \`a des ph\'enom\`enes de semi-continuit\'e 
(lorsque les param\`etres varient). Le but de cet article est 
de fournir un cadre th\'eorique rigoureux \`a cette intuition qui 
nous permettra d'\'enoncer des conjectures pr\'ecises. 
\end{abstract}

\maketitle

\pagestyle{myheadings}

\markboth{\sc C. Bonnaf\'e}{\sc Semi-continuit\'e des cellules de Kazhdan-Lusztig}


Soit $(W,S)$ un groupe de Coxeter (o\`u $S$ est fini) 
et supposons pour simplifier dans cette 
introduction que $S=S_1 \coprod S_2$, o\`u $S_1$ et $S_2$ sont deux 
sous-ensembles non vides de $S$ tels que, si $s_1 \in S_1$ et 
$s_2 \in S_2$, alors $s_1$ et $s_2$ ne sont pas conjugu\'es 
dans $W$. Notons $\ell : W \to \NM$ 
la fonction longueur, $\ell_i : W \to \NM$ la $S_i$-longueur (si 
$w \in W$, $\ell_i(w)$ d\'esigne le nombre d'\'el\'ements de $S_i$ 
apparaissant dans une d\'ecomposition r\'eduite de $w$~: ce nombre 
ne d\'epend pas du choix de la d\'ecomposition r\'eduite). Fixons 
deux entiers naturels non nuls $a$ et $b$ et posons 
$L_{a,b} : W \to \ZM$, $w \mapsto a\ell_1(w)+b\ell_2(w)$. Alors 
$L_{a,b}$ est une fonction de poids (au sens de Lusztig \cite[\S 3.1]{lusztig}) 
et il est donc possible de d\'efinir une partition 
de $W$ en cellules de Kazhdan-Lusztig \cite[chapitre 8]{lusztig}. 
Il est clair que cette partition ne d\'epend que de $b/a$ (et non 
du couple $(a,b)$)~: nous la noterons $\LC_{b/a}(W)$. 
L'exp\'erience sugg\`ere la conjecture 
suivante~:


\begin{quotation}
\noindent{\bf Conjecture 0.} 
{\it Il existe un entier naturel $m$ et 
des nombres rationnels $0 < r_1 < \cdots < r_m$ (ne d\'ependant que de $W$) 
tels que (en notant $r_0=0$ et $r_{m+1}=+\infty$), si $\th$ et 
$\th'$ sont deux nombres rationnels strictement positifs, alors:
\begin{itemize}
\itemth{a} Si $0 \le i \le m$ et si $r_i < \th,\th' < r_{i+1}$, 
alors $\LC_\th(W)=\LC_{\th'}(W)$. 

\itemth{b} Si $1 \le i \le m$ et si $r_{i-1} < \th < r_i < \th' < r_{i+1}$, 
alors $\LC_{r_i}(W)$ est la partition la plus fine de $W$ qui soit 
\`a la fois moins fine que $\LC_\th(W)$ et moins fine que $\LC_{\th'}(W)$. 
\end{itemize}}
\end{quotation}


\noindent{\sc Remarques - } (1) 
On peut \'evidemment \'enoncer des conjectures 
similaires concernant les partitions en cellules \`a droite et 
en cellules bilat\`eres.

\smallskip

(2) Dans le cas o\`u $W$ est fini, l'existence des nombres rationnels 
$0 < r_1 < \cdots < r_m$ v\'erifiant (a) est facile. En revanche, 
m\^eme dans ce cas, 
la propri\'et\'e (b) est encore \`a ce jour une conjecture. 
Une conjecture de Lusztig \cite[\S 13.12]{lusztig} nous autorise \`a 
esp\`erer qu'elle soit encore vraie pour des groupes infinis.\finl

\bigskip

On peut interpr\'eter cette conjecture en termes topologiques~:
les cellules de Kazhdan-Lusztig 
devraient ob\'eir \`a des ph\'enom\`enes de semi-continuit\'e. 
Le but de cet article est de fournir un cadre rigoureux \`a cette 
intuition en g\'en\'eralisant la conjecture 0 
simultan\'ement dans les deux directions suivantes:
\begin{quotation}
\noindent $\bullet$ On peut tr\`es bien imaginer que la partition 
de $S$ comporte plus de deux parties (cela ne se produit pas 
lorsque $W$ est fini et irr\'eductible).

\noindent $\bullet$ On peut aussi s'int\'eresser aux fonctions de 
poids \`a valeurs dans n'importe quel groupe ab\'elien totalement 
ordonn\'e, et dont les valeurs sur les r\'eflexions simples 
ne sont pas n\'ecessairement strictement positives. 
\end{quotation}
Ces deux g\'en\'eralisations nous conduisent \`a d\'efinir des classes 
d'\'equivalences de fonctions de poids (dans le cas facile d\'etaill\'e 
dans l'introduction, cela a \'et\'e fait en introduisant le rapport 
$b/a$ \`a la place du couple $(a,b)$) et \`a d\'efinir une topologie 
sur l'ensemble de ces classes d'\'equivalence. Cela sera fait 
dans un cadre abstrait dans les sections 1, 2 et 3. Dans les sections 4 et 5, 
nous rappelons les d\'efinitions et quelques propri\'et\'es 
des cellules de Kazhdan-Lusztig, en int\'egrant le formalisme des 
sections 1, 2 et 3. 
Nous \'enon\c{c}ons nos conjectures dans la section 6, et nous les 
illustrerons dans la section 7 en d\'etaillant les exemples suivants~:

\bigskip

\noindent{\sc Exemples - } 
(1) {\it Groupes di\'edraux~:} 
Si $|S|=2$ et si $|S_1|=|S_2|=1$, alors la conjecture 0 
est v\'erifi\'ee en prenant $m=1$ et $r_1=1$ 
(voir \cite[\S 8.8]{lusztig}).

\medskip

(2) {\it Type $F_4$~:} Si $(W,S)$ est de type $F_4$ et 
si $|S_1|=|S_2|=2$, alors la conjecture 0 est 
v\'erifi\'ee en prenant $m=3$ et $r_1=1/2$, $r_2=1$ et $r_3=2$ 
(voir \cite[Corollaire 4.8]{geck f4}). 

\medskip

(3) {\it Type $B_n$~:} Supposons que $(W,S)$ est de type 
$B_n$ (avec $n \ge 2$) et $|S_1|=n-1$ et $|S_2|=1$. 
Dans \cite[Conjectures A et B]{bgil}, la conjecture 0 est 
pr\'ecis\'ee~: il devrait suffire de prendre $m=n-1$ et $r_i=i$. 
Cela a \'et\'e v\'erifi\'e pour $n \le 6$ 
(voir \S\ref{section B}). 

\medskip

(4) {\it Type $B(p,q)$~:} Soient $p$ et $q$ deux nombres naturels 
non nuls et supposons que $|S|=p+q$, $S_1=\{s_1,\dots,s_p\}$, 
$S_2=\{t_1,\dots,t_q\}$ et supposons que le graphe de Coxeter 
de $(W,S)$ soit le suivant
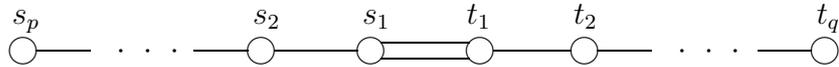
\begin{center}
\begin{picture}(350,30)
\put( 45, 10){\line(1,0){20}}
\put( 40, 10){\circle{10}}
\put( 75,  7){$\cdot$}
\put( 85,  7){$\cdot$}
\put( 95,  7){$\cdot$}
\put(104, 10){\line(1,0){20}}
\put(129, 10){\circle{10}}
\put(134, 10){\line(1,0){31}}
\put(170, 10){\circle{10}}
\put(174,  7){\line(1,0){33}}
\put(174, 13){\line(1,0){33}}
\put(211, 10){\circle{10}}
\put(216, 10){\line(1,0){29}}
\put(250, 10){\circle{10}}
\put(255, 10){\line(1,0){20}}
\put(285,  7){$\cdot$}
\put(295,  7){$\cdot$}
\put(305,  7){$\cdot$}
\put(315, 10){\line(1,0){20}}
\put(340, 10){\circle{10}}
\put( 36, 20){$s_p$}
\put(126, 20){$s_2$}
\put(167, 20){$s_1$}
\put(206, 20){$t_1$}
\put(246, 20){$t_2$}
\put(337, 20){$t_q$}
\end{picture}
\end{center}
Nous dirons qu'il est de type $B(p,q)$. Notons que $B(p,q) \simeq B(q,p)$, 
$B(p,1) \simeq B_{p+1}$, $B(2,2) \simeq F_4$ et 
$B(3,2) \simeq B(2,3) \simeq \tilde{F}_4$. On peut alors se demander si 
les exemples (2) et (3) pr\'ec\'dents ne se g\'en\'eralisent pas ainsi~: 
est-il vrai que la conjecture 0 est valide en prenant $m=p+q-1$ et, 
pour la suite $r_1 < \cdots < r_{p+q-1}$, la suite croissante 
$\DS{\frac{1}{q} < \frac{1}{q-1} < \cdots < \frac{1}{2} < 1 < 2 < \cdots 
< p-1 < p}$~?

\newpage
~
\vskip1cm

\tableofcontents

\newpage

\section{Parties positives d'un r\'eseau\label{section positive}}

\medskip

\begin{quotation}
{\it Fixons dans cette section un r\'eseau $\L$ et notons 
$V=\RM \otimes_\ZM \L$.}
\end{quotation}

\medskip

Le but de cette section est d'\'etudier l'ensemble des parties positives 
(voir \S\ref{sous positive} pour la d\'efinition) de $\L$. Nous munirons 
cet ensemble d'une topologie dans la section suivante. 

\def\tors{{\mathrm{tor}}}

\bigskip

\subsection{D\'efinitions, pr\'eliminaires\label{sous positive}} 
Une partie $X$ de $\L$ est dite {\it positive} si les conditions 
suivantes sont satisfaites~:
\begin{itemize}\itemindent1cm
\itemth{P1} $\L=X \cup (-X)$.

\itemth{P2} $X+X \subset X$.

\itemth{P3} $X \cap (-X)$ est un sous-groupe de $\L$.
\end{itemize}
Nous noterons $\positif(\L)$ l'ensemble des parties 
positives de $\L$. Donnons quelques exemples. 
Pour commencer, notons que $\L \in \positif(\L)$. 
Soit $\G$ un groupe totalement ordonn\'e et soit $\ph : \L \to \G$ 
un morphisme de groupes. Posons 
$$\Pos(\ph)=\{\l \in \L~|~\ph(\l) \ge 0\}$$
$$\Pos^+(\ph)=\{\l \in \L~|~\ph(\l) > 0\}.\leqno{\text{et}}$$
Alors il est clair que 
\equat\label{noyau pos}
\Ker \ph = \Pos(\ph) \cap \Pos(-\ph) = \Pos(\ph) \cap -\Pos(\ph)
\endequat
et que 
\equat\label{positif exemple}
\text{\it $\Pos(\ph)$ est une partie positive de $\L$}.
\endequat

\bigskip

\begin{lem}\label{proprietes positives}
Soit $X$ une partie positive de $\L$. Alors~:
\begin{itemize}
\itemth{a} $-X \in \positif(\L)$. 

\itemth{b} $0 \in X$.

\itemth{c} Si $\l \in \L$ et si $r \in \ZM_{>0}$ est tel que 
$r\l \in X$. Alors $\l \in X$.

\itemth{d} $\L/(X \cap (-X))$ est sans torsion. 
\end{itemize}
\end{lem}

\bigskip

\begin{proof}
(a) est imm\'ediat. 
(b) d\'ecoule de la propri\'et\'e (P1) des parties positives. 
(d) d\'ecoule de (c). Il nous reste \`a montrer (c). 
Soient $\l \in \L$ et $r \in \ZM_{>0}$ tels que $r\l \in X$. 
Si $\l \not\in X$, alors $-\l \in X$ 
d'apr\`es la propri\'et\'e (P1). 
D'o\`u $\l=r\l + (r-1)(-\l) \in X$ d'apr\`es (P2), ce qui est contraire 
\`a l'hypoth\`ese. Donc $\l \in X$. 
\end{proof}

\bigskip

%

Nous allons montrer une forme de r\'eciproque facile \`a la propri\'et\'e 
\ref{positif exemple}. Soit $X \in \positif(\L)$. Notons 
$\can_X : \L \to \L/(X \cap(-X))$ le morphisme canonique. Si $\g$ et $\g'$ 
appartiennent \`a $\L/(X \cap(-X))$, nous \'ecrirons $\g \le_X \g'$ s'il existe 
un repr\'esentant de $\g' - \g$ appartenant \`a $X$. Il est facile 
de v\'erifier que
\equat\label{relation independante}
\text{\it $\g \le_X \g'$ si et seulement si tout repr\'esentant 
de $\g'-\g$ appartient \`a $X$.}
\endequat
On d\'eduit alors facilement des propri\'et\'es (P1), (P2) et (P3) 
des parties positives que 
\equat\label{ordre X}
\text{\it $(\L/(X \cap (-X)),\le_X)$ est un groupe ab\'elien 
totalement ordonn\'e}
\endequat
et que 
\equat\label{X pos}
X=\Pos(\can_X).
\endequat

\bigskip

\subsection{Cons\'equences du th\'eor\`eme de Hahn-Banach} 
Si $X$ est une partie positive de $\L$, on pose 
$X^+=X \setminus (-X)$. On a alors 
\equat\label{disjonction}
\L = X \dotcup (-X^+) = X^+ \dotcup (-X) 
= X^+ \dotcup (X \cap (-X)) \dotcup (-X^+),
\endequat
o\`u $\dot{\cup}$ d\'esigne l'union disjointe. 
De plus, si $\ph : \L \to \G$ est un morphisme de groupes 
ab\'eliens et si $\G$ est un groupe totalement ordonn\'e, alors 
\equat\label{pos plus}
\Pos^+(\ph) = \Pos(\ph)^+.
\endequat
Si $\ph$ est une forme lin\'eaire sur $V$, nous noterons abusivement 
$\Pos(\ph)$ et $\Pos^+(\ph)$ les parties $\Pos(\ph|_\L)$ et 
$\Pos^+(\ph|_\L)$ de $\L$. 

\bigskip

\begin{lem}\label{equivalent separant}
Soit $X$ une partie positive propre de $\L$, soit $\G$ un groupe ab\'elien 
totalement ordonn\'e archim\'edien 
et soit $\ph : \L \to \G$ un morphisme de groupes. 
Alors les conditions suivantes sont \'equivalentes~:
\begin{itemize}
\itemth{1} $X \subseteq \Pos(\ph)$;

\itemth{2} $X^+ \subseteq \Pos(\ph)$;

\itemth{3} $\Pos^+(\ph) \subseteq X^+$;

\itemth{4} $\Pos^+(\ph) \subseteq X$.
\end{itemize}
\end{lem}

\begin{proof}
Il est clair que (1) implique (2) et que (3) implique (4). 

\medskip

Montrons que (2) implique (3). Supposons 
donc que (2) est v\'erifi\'ee. Soit $\l \in \L$ tel que 
$\ph(\l) > 0$ et supposons que $\l \not\in X^+$. 
Alors, d'apr\`es \ref{disjonction}, $\l \in -X$. 
Or, si $\mu \in \L$, il existe $k \in \ZM_{> 0}$ tel que 
$\ph(\mu - k \l)=\ph(\mu)-k\ph(\l) < 0$ (car $\G$ est archim\'edien). Donc $\mu-k\l \not\in X^+$ 
d'apr\`es (2). Donc $\mu - k \l \in -X$ d'apr\`es \ref{disjonction}. 
Donc, $\mu=(\mu -k\l) + k\l \in (-X)$. Donc $\L \subseteq -X$, 
ce qui est contraire \`a l'hypoth\`ese.

\medskip

Montrons que (4) implique (1). 
Supposons donc que $\Pos^+(\ph) \subseteq X$. 
En prenant le compl\'ementaire dans $\L$, 
on obtient $(-X^+)=(-X)^+ \subseteq \Pos(-\ph)$ et donc, 
puisque (2) implique (3), on a $\Pos^+(-\ph) \subseteq -X^+$. 
En reprenant le compl\'ementaire dans $\L$, on 
obtient $X \subseteq \Pos(\ph)$. 
\end{proof}

\bigskip

Nous aurons aussi besoin du lemme suivant~:

\bigskip

\begin{lem}\label{solutions}
Soient $\l_1$,\dots, $\l_n$ des \'el\'ements de $\L$ et supposons 
trouv\'e un $n$-uplet $t_1$,\dots, $t_n$ de nombres r\'eels 
{\bfit strictement} positifs tels que $t_1 \l_1 + \cdots + t_n \l_n=0$. 
Alors il existe des entiers naturels non nuls $r_1$,\dots, $r_n$ tels que 
$r_1 \l_1 + \cdots + r_n \l_n = 0$.
\end{lem}

\bigskip

\begin{proof}
Notons $\SC$ l'ensemble des $n$-uplets $(u_1,\dots,u_n)$ 
de nombres r\'eels qui satisfont 
$$\begin{cases}
u_1+\cdots + u_n = 1, &\\
u_1 \l_1 + \cdots + u_n \l_n = 0.&\\
\end{cases}$$
\'Ecrit dans une base de $\L$, ceci est un syst\`eme lin\'eaire 
d'\'equations \`a coefficients dans $\QM$. Le proc\'ed\'e d'\'elimination de 
Gauss montre que l'existence d'une solution {\it r\'eelle} implique 
l'existence d'une solution {\it rationnelle} 
$t^\circ = (t_1^\circ,\dots, t_n^\circ)$ 
et l'existence de vecteurs $v_1$,\dots, $v_r \in \QM^n$ tels que 
$$\SC=\{t^\circ + x_1 v_1 + \cdots +x_r v_r~|~(x_1,\dots,x_r) \in \RM^r\}.$$
En particulier, il existe $x_1$,\dots, $x_r \in \RM$ tels que 
$$(t_1,\dots,t_n) = t^\circ + x_1 v_1 + \cdots +x_r v_r.$$
Puisque $t_1$,\dots, $t_n$ sont strictement positifs, 
il existe $x_1'$,\dots, $x_r'$ dans $\QM$ tels que les 
coordonn\'ees de $t^\circ + x_1' v_1 + \cdots +x_r' v_r$ soient 
strictement positives. Posons alors 
$(u_1,\dots,u_n)=t^\circ + x_1' v_1 + \cdots +x_r' v_r$. 
On a donc $u_i \in \QM_{> 0}$ pour tout $i$ et 
$$u_1 \l_1 + \cdots + u_n \l_n = 0.$$
Quitte \`a multiplier par le produits des d\'enominateurs des $u_i$, 
on a trouv\'e $r_1$,\dots, $r_n \in \ZM_{> 0}$ tels que 
$$r_1 \l_1 + \cdots + r_n \l_n = 0,$$
comme attendu.
\end{proof}

\bigskip

\begin{theo}\label{hahn banach}
Soit $X$ une partie positive de $\L$ diff\'erente de $\L$. Alors~:
\begin{itemize}
\itemth{a} Il existe une forme lin\'eaire $\ph$ sur $V$ telle que 
$X \subseteq \Pos(\ph)$. 

\itemth{b} Si $\ph$ et $\ph'$ sont deux formes lin\'eaires 
sur $V$ telles que $X \subseteq \Pos(\ph) \cap \Pos(\ph')$, 
alors il existe $\kappa \in \RM_{> 0}$ 
tel que $\ph'=\kappa \ph$.
\end{itemize}
\end{theo}

\bigskip

\begin{proof}
Notons $\CC^+$ l'enveloppe convexe de $X^+$. 
Notons que $X^+$ (et donc $\CC^+$) 
est non vide car $X \neq \L$. Nous allons commencer 
par montrer que $0 \not\in \CC^+$. Supposons donc que $0 \in \CC^+$. 
Il existe donc $\l_1$,\dots, $\l_n$ dans $X^+$ et $t_1$,\dots, $t_n$ 
dans $\RM_{> 0}$ tels que 
$$\begin{cases}
t_1+\cdots + t_n = 1, &\\
t_1 \l_1 + \cdots + t_n \l_n = 0.&\\
\end{cases}$$
D'apr\`es le lemme \ref{solutions}, il existe $r_1$,\dots, $r_n \in \ZM_{> 0}$ 
tels que $m_1 \l_1 = - (m_2 \l_2 + \cdots + m_n \l_n)$. Donc 
$m_1\l_1 \in X \cap -X$ 
(voir la propri\'et\'e (P2)). Donc, d'apr\`es le lemme 
\ref{proprietes positives} (a) et (c), on a $\l_1 \in X \cap -X$, 
ce qui est impossible car $\l_1 \in X^+ = X \setminus (-X)$. 
Cela montre donc que 
$$0 \not\in \CC^+.$$
L'ensemble $\CC^+$ \'etant convexe, il d\'ecoule du th\'eor\`eme de 
Hahn-Banach qu'il existe une forme lin\'eaire non nulle $\ph$ 
sur $V$ telle que
$$\CC^+ \subseteq \{\l \in V~|~\ph(\l) \ge 0\}.$$
En particulier, 
$$X^+ \subseteq \CC^+ \cap \L \subseteq \Pos(\ph).\leqno{(*)}$$
D'apr\`es le lemme \ref{equivalent separant}, et puisque 
$\RM$ est archim\'edien, on a bien $X \subseteq \Pos(\ph)$. 
Cela montre (a). 

\medskip

Soit $\ph'$ une autre forme lin\'eaire telle que $X \subseteq \Pos(\ph')$. 
Posons $U=\{\l \in V~|~\ph(\l) > 0$ et $\ph'(\l) < 0\}$. 
Alors $U$ est un ouvert de $V$. Fixons une $\ZM$-base $(e_1,\dots,e_d)$ 
de $\L$. Si $U \neq \vide$, alors il existe 
$\l \in \QM \otimes_\ZM \L$ et $\e \in \QM_{> 0}$ tels que 
$\l$, $\l + \e e_1$,\dots, $\l+\e e_d$ appartiennent \`a $U$. 
Quitte \`a multiplier par le produit des d\'enominateurs de 
$\e$ et des coordonn\'ees de $\l$ dans la base $(e_1,\dots,e_d)$, 
on peut supposer que $\l \in \L$ et $\e \in \ZM_{> 0}$. Mais 
il est clair que $U \cap X^+=\vide$ et $U \cap (-X^+) = \vide$. 
Donc $U \cap \L$ est contenu dans $X \cap (-X)$. Ce dernier 
\'etant un sous-groupe, on en d\'eduit que $\e e_i \in X \cap (-X)$ 
pour tout $i$. D'o\`u, d'apr\`es le lemme \ref{proprietes positives} (c), 
$e_i \in X \cap (-X)$ pour tout $i$. Donc $X=\L$, ce qui 
est contraire \`a l'hypoth\`ese. On en d\'eduit que $U$ est 
vide, c'est-\`a-dire qu'il existe $\kappa \in \RM_{> 0}$ tel que 
$\ph'=\kappa\ph$. D'o\`u (b).
\end{proof}

\bigskip

Si $\ph$ est une forme lin\'eaire sur $V$, nous noterons 
$\phba$ sa classe dans $V^*/\RM_{>0}$. Nous noterons $p : V^* \longto V^*/\RM_{>0}$
la projection canonique. L'application $\Pos : V^* \to \positif(\L)$ 
se factorise \`a travers $p$ en une application $\Posbar : 
V^*/\RM_{>0} \to \positif(\L)$ 
rendant le diagramme
$$\diagram
V^* \ddrrto^{\DS{\Pos}} \ddto_{\DS{p}}&& \\
&&\\
V^*\!/\RM_{>0} \rrto^{\DS{\Posbar}} && \positif(\L)
\enddiagram$$
commutatif. D'autre part, 
si $X \in \positif(\L)$, nous noterons $\pi(X)$ l'unique 
\'el\'ement $\phba \in V^*/\RM_{>0}$ tel que $X \subseteq \Posbar(\phba)$ 
(voir le th\'eor\`eme \ref{hahn banach}). On a donc d\'efini 
deux applications
$$\diagram
V^*/\RM_{>0}\rrto^{\DS{\Posbar}} && 
\positif(\L) \rrto^{\DS{\pi}} && V^*\!/\RM_{>0} 
\enddiagram$$
et le th\'eor\`eme \ref{hahn banach} (b) montre que
\equat\label{pi pos}
\pi \circ \Posbar = \Id_{V^*\!/\RM_{>0}}.
\endequat
Donc $\pi$ est surjective et $\Posbar$ est injective. 
En revanche, ni $\pi$, ni $\Posbar$ ne sont des bijections (sauf 
si $\L$ est de rang $1$). Nous allons d\'ecrire les fibres de $\pi$~:

\bigskip 

\begin{prop}\label{fibres pi}
Soit $\ph$ une forme lin\'eaire non nulle sur $V$. Alors 
l'application 
$$\fonction{i_{\phba}}{\positif(\Ker \ph|_\L)}{\pi^{-1}(\phba)}{X}{X 
\cup \Pos^+(\ph)}$$
est bien d\'efinie et bijective. Sa r\'eciproque est l'application 
$$\fonctio{\pi^{-1}(\phba)}{\positif(\Ker \ph|_\L)}{Y}{Y \cap \Ker \ph|\l.}$$
\end{prop}

\noindent{\sc Remarque - } 
Il est facile de voir que $\pi^{-1}(\bar{0})=\{\L\}$. Nous 
pouvons aussi d\'efinir une application 
$i_{\bar{0}} : \positif(\L) \to \positif(\L)$ 
par la m\^eme formule que dans la proposition 
\ref{fibres pi}~: alors $i_{\bar{0}}$ est tout simplement 
l'application identit\'e mais on a dans ce cas-l\`a 
$\pi^{-1}(\bar{0})\neq i_{\bar{0}}(\positif(\L))$.\finl

\bigskip

\begin{proof}
Montrons tout d'abord que l'application $i_\phba$ est bien d\'efinie. 
Soit $X \in \positif(\Ker \ph|\L)$. Posons $Y=X \cup \Pos^+(\ph)$. 
Montrons que $Y$ est une partie positive de $\L$. Avant cela, 
notons que 
$$Y \subseteq \Pos(\ph).\leqno{(*)}$$

(1) Si $\l \in \L$, deux cas se pr\'esentent. Si $\ph(\l) \neq 0$, 
alors $\l \in \Pos^+(\ph) \cup -\Pos^+(\ph) \subseteq Y \cup (-Y)$. 
Si $\ph(\l)=0$, alors $\l \in \Ker \ph|_\L$, donc 
$\l \in X \cup (-X) \subseteq Y \cup (-Y)$ car $X$ est 
une partie positive de $\Ker \ph|_\L$. Donc $\L=Y \cup (-Y)$. 

\smallskip

(2) Soient $\l$, $\mu \in Y$. Montrons que $\l+\mu \in Y$. 
Si $\ph(\l+\mu) > 0$, alors $\l+\mu \in \Pos^+(\ph) \subseteq Y$. 
Si $\ph(\l+\mu)=0$, alors il r\'esulte de $(*)$ que 
$\ph(\l)=\ph(\mu)=0$, donc $\l$, $\mu \in \Ker \ph|_\L$. 
En particulier, $\l$, $\mu \in X$ et donc 
$\l+\mu \in X+X \subseteq X \subseteq Y$. 
Donc $Y + Y \subseteq Y$.

\smallskip

(3) On a $Y \cap (-Y)=X \cap (-X)$, donc $Y \cap (-Y)$ est un sous-groupe 
de $\L$. 

\medskip

Les points (1), (2) et (3) ci-dessus montrent que $Y$ est une partie 
positive de $\L$. L'inclusion $(*)$ montre que $\pi(Y)=\phba$, 
c'est-\`a-dire que $Y \in \pi^{-1}(\phba)$. Donc l'application 
$i_\phba$ est bien d\'efinie.

\medskip

Elle est injective car, si $X \in \positif(\Ker \ph|_\L)$, alors 
$X \cap \Pos^+(\ph) = \vide$. Montrons maintenant 
qu'elle est surjective. Soit $Y \in \pi^{-1}(\phba)$. 
Posons $X=Y \cap \Ker \ph|_\L$. Alors $X$ est une partie 
positive de $\Ker \ph|_\L$ d'apr\`es le corollaire 
\ref{pos inclusion}. Posons $Y'=X \cup \Pos^+(\ph)$. 
Il nous reste \`a montrer que $Y=Y'$. 

Tout d'abord, $\Pos^+(\ph) \subseteq Y$ car $\pi(Y)=\phba$ 
par hypoth\`ese 
et $X \subseteq Y$. Donc $Y' \subseteq Y$. R\'eciproquement, 
si $\l \in Y$, deux cas se pr\'esentent. Si $\ph(\l) > 0$, 
alors $\l \in \Pos^+(\ph) \subseteq Y'$. Si $\ph(\l)=0$, 
alors $\l \in Y \cap \Ker \ph|_\L=X \subseteq Y'$. 
Dans tous les cas, $\l \in Y'$.
\end{proof}

\bigskip

\example{maximal} 
Si $\ph$ est une forme lin\'eaire non nulle sur $V$, alors 
$\Pos(\ph) = i_{\phba}(\Ker \ph|_\L)$.\finl

\bigskip

Nous pouvons maintenant classifier les parties positives de 
$\L$ en termes de formes lin\'eaires. Notons $\FC(\L)$ 
l'ensemble des suites finies $(\ph_1,\dots,\ph_r)$ 
telles que (en posant $\ph_0=0$), pour tout $i \in \{1,2,\dots,r\}$, 
$\ph_i$ soit une forme lin\'eaire non nulle sur 
$\RM \otimes_\ZM (\L \cap \Ker \ph_{i-1})$. 
Par convention, nous supposerons que la suite vide, not\'ee $\vide$, 
appartient \`a $\FC(\L)$. 

Posons $d = \dim V$. 
Notons que, si $(\ph_1,\dots,\ph_r) \in \FC(\L)$, alors $r \le d$. 
Nous d\'efinissons donc l'action suivante de $(\RM_{>0})^d$ sur $\FC(\L)$~: 
si $(\kappa_1,\dots,\kappa_d) (\RM_{>0})^d$ et si 
$(\ph_1,\dots,\ph_r) \in \FC(\L)$, on pose 
$$(\kappa_1,\dots,\kappa_d)\cdot (\ph_1,\dots,\ph_r) = 
(\kappa_1 \ph_1,\dots,\kappa_r \ph_r).$$
Munissons $\RM^r$ de l'ordre lexicographique~: c'est un group ab\'elien 
totalement ordonn\'e et $(\ph_1,\dots,\ph_r) : \L \to \RM^r$ 
est un morphisme de groupes. Donc $\Pos(\ph_1,\dots,\ph_r)$ est 
bien d\'efini et appartient \`a $\positif(\L)$. En fait, toute 
partie positive de $\L$ peut \^etre retrouv\'ee ainsi~:

\bigskip

\begin{prop}\label{surjection pos}
L'application 
$$\fonctio{\FC(\L)}{\positif(\L)}{\phb}{\Pos(\phb)}$$
est bien d\'efinie et induit une bijection 
$\FC(\L)/(\RM_{> 0})^d \stackrel{\sim}{\longrightarrow} \positif(\L)$. 
\end{prop}

\noindent{\sc Remarque - } Dans cette proposition, on a pos\'e 
par convention $\Pos(\vide)=\L$.\finl

\bigskip

\begin{proof}
Cela r\'esulte imm\'ediatement d'un raisonnement par r\'ecurrence 
sur le rang de $\L$ en utilisant le th\'eor\`eme \ref{hahn banach} 
et la proposition \ref{fibres pi}.
\end{proof}

\bigskip

Si $\phb \in \FC(\L)$, nous noterons $\bar{\phb}$ sa classe 
dans $\FC(\L)/(\RM_{>0})^d$ et nous poserons 
$\Posbar(\bar{\phb})=\Pos(\phb)$. Comme corollaire de la proposition 
\ref{surjection pos}, nous obtenons une classification des ordres 
totaux sur $\L$ (compatibles avec la structure de groupe). 
En fait, se donner un ordre 
total sur $\L$ est \'equivalent \`a se donner une partie 
positive $X$ de $\L$ telle que $X \cap (-X)=0$. Notons 
$\FC_0(\L)$ l'ensemble des \'el\'ements $(\ph_1,\dots,\ph_r) \in \FC(\L)$ 
tels que $\L \cap \Ker \ph_r = 0$. 

\bigskip

\begin{coro}\label{ordres totaux}
L'application d\'ecrite dans la proposition \ref{surjection pos} induit 
une bijection entre l'ensemble $\FC_0(\L)/(\RM_{> 0})^d$ et l'ensemble 
des ordres totaux sur $\L$ compatibles avec la structure 
de groupe. 
\end{coro}

\bigskip

\subsection{Fonctorialit\'e} 
Soit $\s : \L' \to \L$ un morphisme de groupes ab\'eliens. 
Le r\'esultat suivant est facile~:

\bigskip

\begin{lem}\label{sigma}
Si $X$ est une partie positive de $\L$, 
alors $\s^{-1}(X)$ est une partie positive de $\L'$.
\end{lem}

\begin{proof}
Soit $X \in \positif(\L)$. Posons $X' = \s^{-1}(X)$. 
Montrons que $X' \in \positif(\L')$.

\medskip

(1) On a $X' \cup (-X') = \s^{-1}(X \cup (-X))=\s^{-1}(\L)=\L'$.

\smallskip

(2) Si $\l'$ et $\mu'$ sont deux \'el\'ements de $X'$, alors 
$\s(\l')$ et $\s(\mu')$ appartiennent \`a $X$. Donc $\s(\l')+\s(\mu') \in X$. 
En d'autres termes, $\l' + \mu' \in X'$. Donc $X'+X' \subseteq X'$. 

\smallskip

(3) On a $X' \cap (-X') = \s^{-1}(X \cap (-X))$, donc l'ensemble 
$X' \cap (-X')$ est un sous-groupe de $\L'$.
\end{proof}

\bigskip

\begin{coro}\label{pos inclusion}
Si $\L'$ est un sous-groupe de $\L$ et si $X$ est une partie 
positive de $\L$, alors $X \cap \L'$ est une partie positive de $\L'$. 
\end{coro}

\bigskip

Si $\s : \L' \to \L$ est un morphisme de groupes ab\'eliens, nous noterons 
$\s^* : \positif(\L) \to \positif(\L')$, $X \mapsto \s^{-1}(X)$ induite 
par le lemme \ref{sigma}. Si $\t : \L'' \to \L'$ 
est un morphisme de groupe ab\'eliens, il est alors facile 
de v\'erifier que 
\equat\label{composition}
(\t \circ \s)^* = \s^* \circ \t^*.
\endequat

D'autre part, si on note $V' = \RM \otimes_\ZM \L'$, alors 
$\s$ induit une application $\RM$-lin\'eaire 
$\s_\RM : V' \longto V$ dont nous noterons 
$\lexp{t}{\s_\RM} : V^* \longto V^{\prime *}$ l'application duale 
et $\lexp{t}{\sigba}_\RM : V^*/\RM_{>0} \longto V^{\prime *}/\RM_{> 0}$ 
l'application (continue) induite. 
Il est alors facile de v\'erifier que le diagramme
\equat\label{diag sigma}
\diagram
V^*/\RM_{>0} \rrto^{\DS{\Posbar}} \ddto^{\DS{\lexp{t}{\sigba}_\RM}} && 
\positif(\L) \rrto^{\DS{\pi}} \ddto^{\DS{\s^*}} && 
V^*/\RM_{>0} \ddto^{\DS{\lexp{t}{\sigba}_\RM}} \\
&&&&\\
V^{\prime *}/\RM_{>0} \rrto^{\DS{\Posbar'}} && 
\positif(\L') \rrto^{\DS{\pi'}} && 
V^{\prime *}/\RM_{>0}
\enddiagram
\endequat
est commutatif (o\`u $\Posbar'$ et $\pi'$ sont les analogues de $\Posbar$ et $\pi$ 
pour le r\'eseau $\L'$).

\bigskip

\section{Topologie sur $\positif(\L)$}

\medskip

Dans cette section, nous allons d\'efinir sur $\positif(\L)$ une topologie 
et en \'etudier les propri\'et\'es. Nous allons notamment montrer que la 
plupart des applications introduites dans la section pr\'ec\'edente 
($\Posbar$, $\pi$, $i_{\phba}$,...) sont continues. 

\bigskip

\subsection{D\'efinition}
Si $E$ est une partie de $\L$, 
nous poserons 
$$\UC(E)=\{X \in \positif(\L)~|~X \cap E=\vide\}.$$
Si $\l_1$,\dots, $\l_n$ sont des \'el\'ements de $\L$, nous 
noterons pour simplifier $\UC(\l_1,\dots,\l_n)$ l'ensemble 
$\UC(\{\l_1,\dots,\l_n\})$. Si cela est n\'ecessaire, nous noterons ces ensembles 
$\UC_\L(E)$ ou $\UC_\L(\l_1,\dots,\l_n)$. On a alors
\equat\label{u inter}
\UC(E)=\bigcap_{\l \in E} \UC(\l).
\endequat
Notons que 
\equat\label{u vide}
\UC(\vide)=\positif(\L)\qquad\text{et}\qquad \UC(\L)=\{\vide\}.
\endequat
D'autre part, si $(E_i)_{i \in I}$ est une famille de parties de $\L$, alors 
\equat\label{u intersection}
\bigcap_{i \in I} \UC(E_i) = \UC\bigl(\bigcup_{i \in I} E_i\bigr).
\endequat
Une partie $\UC$ de $\positif(\L)$ sera dite {\it ouverte} si, 
pour tout $X \in \UC$, il existe une partie {\bf finie} $E$ de $\L$ telle 
que $X \in \UC(E)$ et $\UC(E) \subset \UC$. L'\'egalit\'e 
\ref{u intersection} montre que cela d\'efinit bien une topologie 
sur $\positif(\L)$. 

\bigskip

\bigskip

\begin{prop}\label{connexite}
Si $\UC$ est un ouvert de $\positif(\L)$ contenant $\L$, alors 
$\UC=\positif(\L)$. En particulier, $\positif(\L)$ est connexe. 
Si $\L \neq 0$, alors il n'est pas s\'epar\'e. 
\end{prop}

\bigskip

\begin{proof}
Par d\'efinition, il existe une partie finie $E$ de $\L \setminus \L$ 
telle que $\UC(E) \subseteq \UC$. Mais on a forc\'ement $E=\vide$, 
donc $\UC=\positif(\L)$ d'apr\`es \ref{u vide}. D'o\`u le 
r\'esultat. 

Le fait que $\positif(\L)$ 
n'est pas s\'epar\'e (lorsque $\L \neq 0$) en d\'ecoule~: 
le point $\L$ de $\positif(\L)$ ne peut \^etre 
s\'epar\'e d'aucun autre.
\end{proof}

\bigskip

\example{Z} 
L'espace topologique $\positif(\ZM)$ n'a que trois points~: 
$\ZM$, $\ZM_{\ge 0}$ et $\ZM_{\le 0}$. Sur ces trois points, 
seul $\ZM$ est un point ferm\'e et $\ZM_{\ge 0}$ et $\ZM_{\le 0}$ 
sont des points ouverts (en effet, $\{\ZM_{\ge 0}\} = \UC(-1)$ 
et $\{\ZM_{\le 0}\} = \UC(1)$).\finl

\bigskip

Il est clair que la topologie sur $\positif(\L)$ d\'efinie ci-dessus 
est la topologie induite par une topologie sur l'ensemble des parties 
de $\L$ (d\'efinie de fa\c{c}on analogue)~: cette derni\`ere est tr\`es 
grossi\`ere mais sa restriction \`a $\positif(\L)$ est plus int\'eressante. 

\medskip

Nous aurons besoin de la propri\'et\'e suivante des ensembles $\UC(E)$~:

\bigskip

\begin{lem}\label{uce vide}
Soit $E$ une partie {\bfit finie} de $\L$. Alors les assertions suivantes 
sont \'equivalentes~:
\begin{itemize}
\itemth{1} $\UC(E)=\vide$.

\itemth{2} Il existe $n \ge 1$, $\l_1$,\dots, $\l_n \in E$ et 
$r_1$,\dots, $r_n \in \ZM_{>0}$ tels que $\DS{\sum_{i=1}^n r_i\l_i=0}$. 

\itemth{3} Il n'existe pas de forme lin\'eaire $\ph$ sur $V^*$ 
telle que $\ph(E) \subset \RM_{>0}$.
\end{itemize}
\end{lem}

\bigskip

\begin{proof}
S'il existe une forme lin\'eaire $\ph$ sur $V^*$ 
telle que $\ph(E) \subseteq \RM_{>0}$, alors $\Pos(-\ph) \in \UC(E)$, et 
donc $\UC(E) \neq \vide$. Donc (1) $\imp$ (3). 

\medskip

Supposons trouv\'es $\l_1$,\dots, $\l_n \in E$ et $r_1$,\dots, $r_n \in \ZM_{>0}$ 
tels que $r_1\l_1 + \cdots + r_n \l_n = 0$. Alors, si $X \in \UC(E)$, 
on a $-\l_2$,\dots, $-\l_n \in X$. Mais 
$r_1\l_1 = -r_2\l_2 - \cdots - r_n \l_n \in X$, 
donc $\l_1 \in X$, ce qui contredit l'hypoth\`ese. Donc (2) $\imp$ (1).

\medskip

Il nous reste \`a montrer que (3) $\imp$ (2). Supposons que (2) n'est pas vraie. 
Nous allons montrer qu'alors (3) n'est pas vraie en raisonnant par r\'ecurrence sur
la dimension de $V$ (c'est-\`a-dire le rang de $\L$). 
Posons 
$$\CC=\{t_1\l_1+\cdots t_n \l_n~|~n \ge 1,
~\l_1,\dots, \l_n \in \L,~t_1,\dots,t_n \in \RM_{>0}\}.$$
Alors $\CC$ est une partie convexe de $V$ contenant $E$ et, d'apr\`es 
le lemme \ref{solutions} et le fait que (2) ne soit pas vraie, on a 
$0 \not\in \CC$. Par cons\'equent, il r\'esulte du th\'eor\`eme de Hahn-Banach 
qu'il existe une forme lin\'eaire non nulle $\ph$ telle que 
$\ph(\CC)\subseteq \RM_{\ge 0}$. Posons $\L'=(\Ker \ph)\cap \L$ et 
$E' = E \cap \L$. Alors l'assertion (2) pour $E'$ n'est pas vraie elle aussi donc, 
par hypoth\`ese de r\'ecurrence, il existe un forme lin\'eaire $\psi$ sur 
$V'=\RM \otimes_\ZM \L' \subseteq V$ telle que $\psi(E') \subset \RM_{>0}$. 
Soit $\psit$ une extension de $\psi$ \`a $V$. 
Puisque $\ph(E \setminus E') \subset \RM_{>0}$, il existe $\e > 0$ tel que 
$\ph(\l) + \e \psit(\l) > 0$ pour tout $\l \in E\setminus E'$. 
Mais on a aussi, si $\l \in E'$, $\ph(\l)+\e\psit(\l) = \e \psit(\l) > 0$. 
Donc $(\ph + \e \psit)(E) \subset \RM_{>0}$.
\end{proof}

\bigskip

Si $\EC$ est une partie de $\positif(\L)$, nous noterons 
$\overline{\EC}$ son adh\'erence dans $\positif(\L)$.

\bigskip

\begin{coro}\label{adherence u}
Soit $E$ une partie {\bfit finie} de $\L$ telle que $\UC(E) \neq \vide$. Alors
$$\overline{\UC(E)} = \positif(\L) \setminus \Bigl(\bigcup_{\l \in E} \UC(-\l)\Bigr).$$
\end{coro}

\bigskip

\begin{proof}
Posons 
$$\OC=\bigcup_{\l \in E} \UC(-\l)$$
$$\FC=\positif(\L) \setminus \OC. \leqno{\text{et}}$$
Alors $\FC$ est ferm\'e dans $\positif(\L)$ et contient $\UC(E)$. 
Donc $\overline{\UC(E)} \subseteq \FC$. 

R\'eciproquement, soit $X \in \positif(\L) \setminus \overline{\UC(E)}$. 
Nous devons montrer que 
$$X \in \OC.\leqno{(?)}$$
Puisque $\positif(\L) \setminus \overline{\UC(E)}$ est un ouvert, il existe 
une partie finie $F$ de $\L$ telle que $X \in \UC(F)$ et 
$\UC(F) \subseteq \OC$. En particulier, $\UC(F) \cap \UC(E) = \vide$. 
En d'autres termes, d'apr\`es \ref{u intersection}, on a 
$\UC(E \cup F) = \vide$. 
Donc, d'apr\`es le lemme \ref{uce vide}, il existe $m \ge 0$, $n \ge 0$, 
$\l_1$,\dots, $\l_m \in E$, $\mu_1$,\dots, $\mu_n \in F$, 
$r_1$,\dots, $r_m$, $s_1$,\dots, $s_n \in \ZM_{>0}$ tels que
$$r_1\l_1 + \cdots + r_m \l_m + s_1 \mu_1 + \cdots + s_n \mu_n=0$$
et $m+n \ge 1$. En fait, comme $\UC(E)$ et $\UC(F)$ sont toutes deux non vides, 
il d\'ecoule du lemme \ref{uce vide} que $m$, $n \ge 1$. 

Si $X \not\in \OC$, alors $-\l_i \in X$ pour tout $i$, ce qui implique 
que $s_1 \mu_1 + \cdots + s_n \mu_n \in X$. Cela ne peut se produire que si au moins 
l'un des $\mu_j$ appartient \`a $X$, mais c'est impossible car $F \cap X =\vide$. 
D'o\`u (?).
\end{proof}

\bigskip

\example{viditude} 
Le corollaire \ref{adherence u} n'est pas forc\'ement vrai si $\UC(E)=\vide$. 
En effet, si $\l \in \L \setminus\{0\}$, alors $\UC(\l,-\l)=\vide$ mais, 
du moins lorsque $\dim V \ge 2$, on a $\UC(\l) \cup \UC(-\l) \neq \positif(\L)$.\finl

\bigskip

\subsection{Fonctorialit\'e} 
Si $\L'$ est un autre r\'eseau et si $\s : \L' \to \L$ est un 
morphisme de groupes et si $E'$ est une partie de $\L'$, alors 
\equat\label{sigma continue}
(\s^*)^{-1}\bigl(\UC_{\L'}(E')\bigr)=\UC_\L\bigl(\s(E')\bigr).
\endequat
\begin{proof}[Preuve de \ref{sigma continue}] 
Soit $X$ une partie positive de $\L$. Alors 
$X \in (\s^*)^{-1}\bigl(\UC_{\L'}(E')\bigr)$ 
(resp. $X \in \UC_\L\bigl(\s(E')\bigr)$) 
si et seulement si $\s^{-1}(X) \cap E' = \vide$ 
(resp. $X \cap \s(E') = \vide$). 
Il est alors facile de v\'erifier que ces deux derni\`eres conditions sont 
\'equivalentes.
\end{proof}

Cela implique le r\'esultat suivant~:

\bigskip

\begin{prop}\label{continu}
L'application $\s^* : \positif(\L) \to \positif(\L')$ est continue. 
\end{prop}

\bigskip

%

\subsection{Continuit\'e} 
D'apr\`es la section \ref{section positive}, nous avons \'equip\'e 
l'espace topologique $\positif(\L)$ de 
deux applications $\Posbar : V^*/\RM_{> 0} \to \positif(\L)$ et 
$\pi : \positif(\L) \to V^*/\RM_{>0}$ telles que 
$\pi \circ \Posbar = \Id_{V^*/\RM_{>0}}$. Nous montrerons dans la 
proposition \ref{pi pos continues} que ces applications 
sont continues (lorsque $V^*/\RM_{>0}$ est bien s\^ur muni de 
la topologie quotient) et nous en d\'eduirons quelques autres 
propri\'et\'es topologiques de ces applications. 
Avant cela, introduisons la notation suivante~: si $E$ 
est une partie finie de $\L$, on pose 
$$\VC(E)=\{\phba\in V^*/\RM_{> 0}~|~\forall~\l \in E,~\ph(\l) < 0\}.$$
Si cela est n\'ecessaire, nous noterons $\VC_\L(E)$ l'ensemble $\VC(E)$. 
Alors
$$p^{-1}(\VC(E))=\{\ph \in V^*~|~\forall~\l \in E,~\ph(\l) < 0\}.$$
Donc $p^{-1}(\VC(E))$ est ouvert, donc $\VC(E)$ est ouvert 
dans $V^*/\RM_{>0}$ par d\'efinition de la topologie quotient. 

\bigskip

\begin{prop}\label{pi pos continues}
Les applications $\Posbar$ et $\pi$ sont continues. De plus~:
\begin{itemize}
\itemth{a} $\Posbar$ induit un hom\'eomorphisme sur son image.

\itemth{b} L'image de $\Posbar$ est dense dans $\positif(\L)$.
%
\end{itemize}
\end{prop}

\begin{proof} 
Soit $E$ une partie finie de $\L$. Alors 
\equat\label{pos uc}
\Posbar^{-1}(\UC(E))=\VC(E).
\endequat
Donc $\Posbar^{-1}(\UC(E))$ est un ouvert de $V^*/\RM_{>0}$. Donc 
$\Posbar$ est continue. 

Montrons maintenant que $\pi$ est continue. Nous proc\`ederons 
par \'etapes~:

\medskip

\begin{quotation}
\begin{lem}\label{base ouverts}
Les $\VC(E)$, o\`u $E$ parcourt l'ensemble des parties 
finies de $\L$, forment une base d'ouverts de $V^*/\RM_{>0}$. 
\end{lem}

\begin{proof}[Preuve du lemme \ref{base ouverts}]
Soit $\UC$ un ouvert de $V^*/\RM_{>0}$ et soit $\ph$ une forme lin\'eaire 
sur $V$ telle que $\phba \in \UC$. 
Nous devons montrer qu'il existe une partie finie $E$ de $V$ 
telle que $\phba \in \VC(E)$ et $\VC(E) \subset \UC$. 

Si $\ph=0$, alors $\UC=V^*/\RM_{>0}$ et le r\'esultat est clair. 
Nous supposerons donc que $\ph \neq 0$. Il existe alors $\l_0 \in \L$ 
tel que $\ph(\l_0) > 0$. Quitte \`a remplacer $\ph$ par un multiple 
positif, on peut supposer que $\ph(\l_0)=1$. Notons $\HC_0$ l'hyperplan 
affine $\{\psi \in V^*~|~\psi(\l_0)=1\}$. Alors l'application naturelle 
$\HC_0 \to V^*/\RM_{>0}$ induit un hom\'eomorphisme 
$\nu : \HC_0 \stackrel{\sim}{\longrightarrow} \VC(-\l_0)$. 
De plus, $\ph=\nu^{-1}(\ph) \in \HC_0$. Donc $\ph \in \nu^{-1}(\UC \cap \VC(-\l_0))$. 
Il suffit dons de v\'erifier que les intersections finies de demi-espaces 
ouverts {\it rationnels} (i.e. de la forme $\{\psi \in \HC_0~|~\psi(\l) > n\}$ o\`u 
$n \in \ZM$ et $\l \in \L \setminus \ZM \l_0$) 
forment une base de voisinages de l'espace affine $\HC_0$, ce qui est imm\'ediat.
\end{proof}
\end{quotation}

\medskip

Compte tenu du lemme \ref{base ouverts}, il suffit de montrer que, 
si $E$ est une partie finie de $\L$, alors $\pi^{-1}(\VC(E))$ 
est un ouvert de $\positif(\L)$. De plus, 
$$\VC(E)=\bigcap_{\l \in E} \VC(\l).$$
Par cons\'equent, la continuit\'e de $\pi$ 
d\'ecoulera du lemme suivant:

\medskip

\begin{quotation}
\begin{lem}\label{image inverse}
Si $\l \in \L$, alors $\pi^{-1}(\VC(\l))$ est un ouvert de 
$\positif(\L)$.
\end{lem}

\begin{proof}[Preuve du lemme \ref{image inverse}]
Soit $X \in \pi^{-1}(\VC(\l))$ et soit $\ph=\pi(X)$. 
Par d\'efinition, $\ph(\l) < 0$ et donc $\l \not\in X$. 
Soit $e_1$,\dots, $e_n$ une $\ZM$-base de $\L$. Il existe un entier 
naturel non nul $N$ tel que $\ph(\l \pm \DS{\frac{1}{N}e_i}) < 0$ pour 
tout $i$. Quitte \`a remplacer $\l$ par $N\l$, on peut supposer que 
$\ph(\l \pm e_i) < 0$ pour tout $i$. On pose alors 
$$E=\{\l+e_1,\l-e_1,\dots,\l+e_n,\l-e_n\}.$$
Alors $X \in \UC(E)$ par construction. Il reste \`a montrer 
que $\UC(E) \subseteq \pi^{-1}(\VC(\l))$. Soit $Y \in \UC(E)$ 
et posons $\psi = \pi(Y)$. Supposons de plus que $\psi \not\in \VC(\l)$. 
On a alors $\psi(\l) \ge 0$. D'autre part, $\psi(\l \pm e_i) \le 0$ pour tout $i$. 
Cela montre que $2\psi(\l)=\psi(\l+e_1)+\psi(\l-e_1) \le 0$, et donc 
$\psi(\l)=\psi(\l+e_i)=0$, et donc $\psi(e_i)=0$ pour tout $i$. Donc $\psi$ 
est nulle et donc $Y=\L$, ce qui contredit le fait que 
$Y \in \UC(E)$. Cela montre donc que $\psi \in \VC(\l)$, comme attendu.
\end{proof}
\end{quotation}

\medskip

Puisque $\pi$ et $\Posbar$ sont continues et v\'erifient 
$\pi \circ \Posbar = \Id_{V^*/\RM_{>0}}$, $\Posbar$ induit un 
hom\'eomorphisme sur son image. D'o\`u (a).

\medskip

L'assertion (b) d\'ecoule du lemme suivant (qui est une cons\'equence imm\'ediate 
de l'\'equivalence entre (1) et (3) dans le lemme \ref{uce vide}) et de \ref{pos uc}~:
\begin{quotation}
\begin{lem}\label{non vide}
Soit $E$ une partie finie de $\L$ telle que $\UC(E) \neq \vide$. 
Alors $\VC(E) \neq \vide$. 
\end{lem}

%
\end{quotation}

La preuve de la proposition \ref{pi pos continues} est termin\'ee. 
\end{proof}

\bigskip

Nous allons maintenant \'etudier les propri\'et\'es topologiques des 
applications $i_\phba$. 

\bigskip

\begin{prop}\label{proprietes topologie}
Soit $\ph \in V^*$ et supposons $\ph \neq 0$. Alors~:
\begin{itemize}
\itemth{a} 
$\DS{\bigcap_{\stackrel{\UC {\mathrm{~ouvert~de~}}\positif(\L)}{\Pos(\ph) \in \UC}} 
\UC = i_{\phba}\bigl(\positif(\Ker \ph|_\L)\bigr)} = \pi^{-1}(\phba)$.

\itemth{d} $i_{\phba}$ est continue et induit un hom\'eomorphisme sur son image.
\end{itemize}
\end{prop}

\begin{proof}
(a) Notons $I_\ph$ l'image de $i_\phba$. On a alors, d'apr\`es 
le lemme \ref{equivalent separant}, 
$$I_\ph=\{X \in \positif(\L)~|~\Pos^+(\ph) \subseteq X\}.\leqno{(*)}$$
Si $\UC$ est un ouvert contenant $\Pos(\ph)$, 
alors il existe une partie finie $E$ de $\L \setminus \Pos(\ph)$ 
telle que $\UC(E) \subseteq \UC$. Mais, si $X$ est dans l'image de 
$i_\phba$, alors $X \subseteq \Pos(\ph)$, donc $X \cap E = \vide$. 
Et donc $X \in \UC$, ce qui montre que 
$$I_\ph \subseteq
\bigcap_{\stackrel{\UC {\mathrm{~ouvert~de~}}\positif(\L)}{\Pos(\ph) \in \UC}} 
\UC.$$
Montrons l'inclusion r\'eciproque. Soit $X \in \positif(\L)$ tel 
que $X \not\in I_\ph$. Posons $\psi=\pi(X)$. Alors $\psiba \neq \phba$ 
donc, d'apr\`es la preuve du th\'eor\`eme \ref{hahn banach}, il existe 
$\l \in \L$ tel que $\ph(\l) < 0$ et $\psi(\l) > 0$. On a donc, 
d'apr\`es le lemme \ref{equivalent separant}, $X \not\in \UC(\l)$. 
D'autre part, $\Pos(\ph) \in \UC(\l)$. D'o\`u (a).

\medskip

Montrons (b). 
On note $\pi_\phba : I_\ph \to \positif(\Ker \ph|_\L)$, 
$X \mapsto X \cap \Ker \ph|_\L$. D'apr\`es la 
proposition \ref{fibres pi}, $\pi_\phba$ est la bijection 
r\'eciproque de $i_\phba : \positif(\Ker \ph|_\L) \to I_\ph$. 
Il nous faut donc montrer que $i_\phba$ et $\pi_\phba$ 
sont continues. Si $F$ est une partie finie de 
$\Ker \ph|_\L$, nous noterons $\UC_\phba(F)$ l'analogue de 
l'ensemble $\UC(F)$ d\'efini \`a l'int\'erieur de $\positif(\Ker \ph|_\L)$. 

Soit $E$ une partie finie de $\L$. Nous voulons montrer 
que $i_\phba^{-1}(\UC(E))$ est un ouvert de $\positif(\Ker \ph|_\L)$. 
S'il existe $\l \in E$ tel que 
$\ph(\l) > 0$, alors $\UC(E) \cap I_\ph = \vide$ (voir $(*)$). On peut 
donc supposer que $\ph(\l) \le 0$ pour tout $\l \in E$. 
Il est alors facile de v\'erifier que 
$$i_\phba^{-1}(\UC(E)) = \UC_\phba(E \cap \Ker \ph|_\L).$$
Donc $i_\phba$ est continue. 

Soit $F$ une partie finie de $\Ker \ph|_\L$. Alors 
$$\pi_\phba^{-1}(\UC_\phba(F)) =\UC(F) \cap I_\ph.$$
Donc $\pi_\phba$ est continue. 
\end{proof}

Nous allons r\'esumer dans le th\'eor\`eme suivant 
la plupart des r\'esutats obtenus dans cette sous-section.

\bigskip

\begin{theo}\label{theo:topologie}
Supposons $\L \neq 0$ et soit $\ph \in V^*$, $\ph \neq 0$. 
\begin{itemize}
\itemth{a} $\positif(\L)$ est connexe. Il n'est pas s\'epar\'e si $\L \neq 0$.

\itemth{b} Les applications 
$\pi : \positif(\L) \to V^*/\RM_{>0}$ et $\Posbar : V^*/\RM_{>0} \to \positif(\L)$ 
sont continues et v\'erifient $\pi \circ \Posbar = \Id_{V^*/\RM_{>0}}$.

\itemth{c} $\Posbar$ induit un hom\'eomorphisme sur son image~; cette 
image est dense dans $\positif(\L)$.

\itemth{d} $\pi^{-1}(\phba)$ est l'intersection des voisinages de $\Pos(\ph)$ 
dans $\positif(\L)$.

\itemth{e} $i_\phba$ est un hom\'eomorphisme.
\end{itemize}
\end{theo}

\bigskip

\section{Arrangements d'hyperplans}

\medskip

L'application continue $\Pos : V^* \to \positif(\L)$ 
a une image dense. Nous allons \'etudier ici comment se 
transpose la notion d'arrangement d'hyperplans (et les 
objets attach\'es~: facettes, chambres, support...) 
\`a l'espace topologique $\positif(\L)$ \`a travers $\Pos$. 
Cela nous permettra d'\'enoncer les conjectures sur 
les cellules de Kazhdan-Lusztig sous la forme la plus g\'en\'erale 
possible. 

\bigskip

\subsection{Sous-espaces rationnels} 
Si $E$ est une partie de $\L$, on pose 
$$\LC(E)=\{X \in \positif(\L)~|~E \subset X \cap (-X)\}.$$
Si cela est n\'ecessaire, nous le noterons $\LC_\L(E)$. 
On appelle {\it sous-espace rationnel} de $\positif(\L)$ 
toute partie de $\positif(\L)$ de la forme $\LC(E)$, o\`u 
$E$ est une partie de $\L$. Si $\l \in \L\setminus \{0\}$, 
on notera $\HC_\l$ le sous-espace rationnel $\LC(\{\l\})$~: 
un tel sous-espace rationnel sera appel\'e un {\it hyperplan 
rationnel}. Notons que 
\equat\label{partage espace}
\positif(\L)= \UC(\l) \dotcup \HC_\l \dotcup \UC(-\l).
\endequat
La proposition suivante justifie quelque peu la terminologie~:

\begin{prop}\label{topo espace}
Soit $E$ est une partie de $\L$. Notons $\L(E)$ le sous-r\'eseau 
$\L \cap \sum_{\l \in E} \QM E$ de $\L$ et soit $\s_E : \L \to \L(E)$ 
l'application canonique. Alors~:
\begin{itemize}
\itemth{a} $\LC(E)=\DS{\bigcap_{\l \in E\setminus\{0\}} \HC_\l} = 
\{X \in \positif(\L)~|~\L(E) \subseteq X\}$.

\itemth{b} $\LC(E)$ est ferm\'e dans $\positif(\L)$. 

\itemth{c} $\Pos^{-1}(\LC(E))=\{\ph \in V^*~|~\forall~\l \in E,~\ph(\l)=0\}=E^\perp$.

\itemth{d} $\Posbar(\pi(\LC(E)) \subseteq \LC(E)$. 

\itemth{e} $\Posbar^{-1}(\LC(E))=\pi(\LC(E))$.

\itemth{f} $\LC(E)=\overline{\Pos(\Pos^{-1}(\LC(E)))}$ .

\itemth{g} L'application $\s_E^* : \positif(\L/\L(E)) \to \positif(\L)$ 
a pour image $\LC(E)$ et induit un hom\'eomorphisme 
$\positif(\L/\L(E)) \stackrel{\sim}{\longrightarrow} \LC(E)$. 
\end{itemize}
\end{prop}

\begin{proof}
La premi\`ere \'egalit\'e de (a) est imm\'ediate. La deuxi\`eme 
d\'ecoule de la proposition \ref{proprietes positives} (c). 
(b) d\'ecoule de (a) et de \ref{partage espace}.
(c) est tout aussi clair. 

\medskip

(d) Si $X \in \Posbar(\pi(\LC(E))$, alors il existe 
$Y \in \LC(E)$ tel que $X=\Posbar(\pi(Y))$. Posons $\phba = \pi(Y)$, 
o\`u $\ph \in V^*$. Alors $E \subseteq Y \cap (-Y)$ et 
$Y \subseteq \Pos(\ph)$. Or, $\ph(\l)=0$ si $\l \in Y \cap (-Y)$, 
donc $X=\Pos(\ph) \in \LC(E)$. 

\medskip

(e) D'apr\`es (d), on a $\pi(\L(E)) \subseteq \Posbar^{-1}(\LC(E))$. 
R\'eciproquement, soit $\ph$ un \'el\'ement de $\Pos^{-1}(\LC(E))$, alors 
$\phba=\pi(\Posbar(\phba)) \in \pi(\LC(E))$. D'o\`u (e).

\medskip

(f) Notons $\FC_E=\Pos(\Pos^{-1}(\LC(E)))$. 
On a $\FC_E \subseteq \LC(E)$ donc 
il d\'ecoule du (a) que $\overline{\FC}_E \subseteq \LC(E)$. 
R\'eciproquement, soit $F$ une partie finie de $\L$ telle que $\UC(F) \cap \FC_E = \vide$. 
Nous devons montrer que $\UC(F) \cap \LC(E)=\vide$. Or, le fait que $\UC(F) \cap \FC_E=\vide$ 
est \'equivalent \`a l'assertion suivante (voir (c))~:
$$\forall~\ph \in E^\perp,~\forall~\l \in F,~\ph(\l) \ge 0.$$
Or, si $\ph \in E^\perp$, alors $-\ph \in E^\perp$, ce qui implique que~:
$$\forall~\ph \in E^\perp,~\forall~\l \in F,~\ph(\l) = 0.$$
En d'autres termes, $F \subseteq (E^\perp)^\perp \cap \L = \L(E)$. Mais, 
si $X \in \LC(E)$, alors $\L(E) \subseteq X$ d'apr\`es (a). Donc 
$X \not\in \UC(F)$, comme esp\'er\'e.

\medskip

(g) Le fait que l'image de $\s_E^*$ soit 
$\LC(E)$ d\'ecoule de (a). D'autre part, $\s_E^*$ est continue d'apr\`es la 
proposition \ref{continu}. Notons 
$$\fonction{\g_E}{\LC(E)}{\positif(\L/\L(E))}{X}{X/\L(E).}$$
Alors $\g_E$ est la r\'eciproque de $\s_E^*$. Il ne nous reste 
qu'\`a montrer que $\g_E$ est continue. Soit donc $F$ une partie finie de $\L/\L(E)$ 
et notons $\Fti$ un ensemble de repr\'esentants des \'el\'ements de $F$ dans $\L$.
On a 
\eqna
\g_E^{-1}(\UC_{\L/\L(E)}(F))&=&\{X \in \LC(E)~|~\forall~\l \in F,~\l \not\in X/\L(E)\}\\
&=&\{X \in \positif(\L)~|~\L(E) \subseteq X\text{ et }\forall~\l \in F,~\l \not\in X/\L(E)\}\\
&=&\{X \in \positif(\L)~|~\L(E) \subseteq X\text{ et }\forall~\l \in \Fti,~\l \not\in X\}\\
&=&\{X \in \LC(E)~|~\forall~\l \in \Fti,~\l \not\in X\}\\
&=&\LC(E) \cap \UC_\L(\Fti).
\endeqna
Donc $\g_E^{-1}(\UC_{\L/\L(E)})$ est un ouvert de $\LC(E)$. Cela montre la continuit\'e de 
$\g_E$.
\end{proof}

\bigskip

\subsection{Demi-espaces} 
Soit $\HC$ un hyperplan rationnel de $\positif(\L)$ et soit $\l\in \L\setminus\{0\}$ 
tel que $\HC=\HC_\l$. D'apr\`es \ref{partage espace}, 
l'hyperplan $\HC$ nous d\'efinit une unique relation d'\'equivalence 
$\smile_{\!\HC}$ sur $\positif(\L)$ pour laquelle les classes d'\'equivalence 
sont $\UC(\l)$, $\HC$ et $\UC(-\l)$~: notons que cette relation 
ne d\'epend pas du choix de $\l$. De plus~:

\bigskip

\begin{prop}\label{composantes connexes}
$\HC$ est un ferm\'e de $\positif(\L)$ et $\UC(\l)$ et $\UC(-\l)$ 
sont les composantes connexes de $\positif(\L)\setminus \HC$. De plus
$$\overline{\UC(\l)}=\UC(\l) \cup \HC_\l.$$
\end{prop}

\bigskip

\begin{proof}
La derni\`ere assertion est un cas particulier du corollaire 
\ref{adherence u}. 

\medskip

Montrons pour finir que $\UC(\l)$ est connexe. Soient $\UC$ et $\VC$ deux ouverts 
de $\UC(\l)$ tels que $\UC(\l)=\UC \coprod \VC$. Alors 
$$\Pos^{-1}(\UC(\l))=\Pos^{-1}(\UC) \coprod \Pos^{-1}(\VC).$$
Mais $\Pos^{-1}(\UC(\l))=\{\ph \in V^*~|~\ph(\l) < 0\}$. Donc 
$\Pos^{-1}(\UC(\l))$ est connexe. Puisque $\Pos$ est continue, 
cela implique que $\Pos^{-1}(\UC)=\vide$ ou $\Pos^{-1}(\VC)=\vide$. 
Le lemme \ref{non vide} implique que $\UC=\vide$ ou $\VC=\vide$.
\end{proof}

\bigskip

Si $X \in \positif(\L)$, nous noterons $\DC_\HC(X)$ la classe 
d'\'equivalence de $X$ sous la relation $\smile_{\! \HC}$. Il r\'esulte 
de la proposition \ref{composantes connexes} que $\overline{\DC_\HC(X)}$ est 
une r\'eunion de classes d'\'equivalences pour $\smile_{\! \HC}$.

\bigskip

\subsection{Arrangements} 
Nous travaillerons d\'esormais sous l'hypoth\`ese suivante~:

\bigskip

\begin{quotation}
{\it Fixons maintenant, et ce jusqu'\`a la fin de cette section, 
un ensemble {\it fini} $\AG$ d'hyperplans rationnels de $\positif(\L)$.}
\end{quotation}

\bigskip

Nous allons red\'efinir, dans notre espace $\positif(\L)$, les notions 
de {\it facettes}, {\it chambres} et {\it faces} associ\'ees \`a $\AG$, 
de fa\c{c}on analogue \`a ce qui se fait pour les arrangements 
d'hyperplans dans un espace r\'eel 
\cite[Chapitre V, \S 1]{bourbaki}. Les propri\'et\'es des 
applications $\pi$ et $\Posbar$ \'etablies pr\'ec\'edemment 
permettent facilement de d\'emontrer des r\'esultats analogues 
en copiant presque mot \`a mot les preuves de 
\cite[Chapitre V, \S 1]{bourbaki}. 

Nous d\'efinissons la relation $\smile_\AG$ sur $\positif(\L)$ de la fa\c{c}on suivante~:
si $X$ et $Y$ sont deux \'el\'ements de $\positif(\L)$, nous \'ecrirons 
$X \smile_\AG Y$ si $X \smile_{\!\HC} Y$ pour tout $\HC \in \AG$. 
Nous appellerons {\it facettes} (ou {\it $\AG$-facettes}) les classes 
d'\'equivalence pour la relation $\smile_\AG$. Nous appellerons {\it chambres} (ou
{\it $\AG$-chambres}) les facettes qui ne rencontrent aucun hyperplan de
$\AG$. Si $\FC$ est une facette, nous noterons
$$\LC_{\AG} ( \FC) = \bigcap_{\overset{\HC \in \AG}{\FC \subset \HC}} \HC,$$
avec la convention habituelle que $\LC_{\AG} ( \FC) = \positif (\L)$ si
$\FC$ est une chambre. Nous l'appellerons le {\it support} de $\FC$ et
nous appellerons {\it dimension} de $\FC$ l'entier 
$$\dim \FC = \dim_{\RM} \Pos^{- 1} ( \LC_{\AG} ( \FC)) .$$
De m\^eme nous appellerons {\it codimension} de $\FC$ l'entier
$$\codim \FC=\dim_\RM V - \dim \FC.$$
Avec ces d\'efinitions, une chambre est une facette de codimension $0$.

\bigskip

\begin{prop}\label{facettes}
Soit $\FC$ une facette et soit $X \in \FC$. Alors~:
\begin{itemize}
\itemth{a} $\FC=\DS{\bigcap_{\HC \in \AG}} \DC_\HC(X)$.

\itemth{b} $\overline{\FC}=\DS{\bigcap_{\HC \in \AG}} \overline{\DC_\HC(X)}$.

\itemth{c} $\overline{\FC}$ est la r\'eunion de $\FC$ et de facettes de dimension 
strictement inf\'erieures.

\itemth{d} Si $\FC'$ est une facette telle que $\overline{\FC}=\overline{\FC}'$, 
alors $\FC=\FC'$.
\end{itemize}
\end{prop}

\bigskip

\begin{proof}
(a) est une cons\'equence des d\'efinitions. Montrons (b). Posons 
$$\AG_1=\{\HC \in \AG~|~\FC \subseteq \HC\}$$
$$\AG_2=\AG \setminus \AG_1.$$
Pour tout $\HC \in \AG$, on fixe un \'el\'ement $\l(\HC) \in \L$ tel que 
$\HC=\HC_{\l(\HC)}$~: si de plus $\HC \in \AG_2$, on choisit $\l(\HC)$ de sorte 
que $\FC \subseteq \UC(\l(\HC))$. On pose 
$$E_i=\{\l(\HC)~|~\HC \in \AG_i\}.$$ 
Par cons\'equent, 
$$\FC=\LC \cap \UC(E_i).$$
Puisque $\LC$ est ferm\'e, $\overline{\FC}$ est aussi l'adh\'erence de $\FC$ dans 
$\LC$. En utilisant alors l'hom\'eomorphisme 
$\s_{E_1}^* : \positif(\L/\L(E_1)) \stackrel{\sim}{\longrightarrow} \LC$ 
de la proposition \ref{topo espace} (g), 
on se ram\`ene \`a calculer l'adh\'erence de $\s_{E_1}^{* -1}(\FC)$ 
dans $\positif(\L/\L(E_1))$. Mais 
$$\s_{E_1}^{* -1}(\FC)=
\s_{E_1}^{* -1}(\UC_\L(E_2))=\UC_{\L/\L(E_1)}(\s_{E_1}(E_2)).$$
Or, d'apr\`es le corollaire \ref{adherence u}, on a 
$$\overline{\UC_{\L/\L(E_1)}(\s_{E_1}(E_2))} = 
\positif(\L/\L(E_1)) \setminus 
\bigl(\bigcup_{\l \in E_2} \UC_{\L/\L(E_1)}(\s_{E_1}(-\l))\bigr).$$
Par cons\'equent,
$$\overline{\FC}=\LC \cap \s_{E_1}^*\Bigl(\positif(\L/\L(E_1)) \setminus 
\bigl(\bigcup_{\l \in E_2} \UC_{\L/\L(E_1)}(\s_{E_1}(-\l))\bigr)\Bigr),$$
et donc
$$\overline{\FC}=\LC \cap 
\Bigl(\bigcap_{\l \in E_2} \bigl(\positif(\L)\setminus
\UC_\L(-\l)\bigr) = 
\bigcap_{\l \in E_1 \cup E_2} \overline{\DC_{\HC_\l}(X)},$$
comme attendu.

\medskip

Montrons maintenant (c). D'apr\`es (b), $\overline{\FC}$ est bien une r\'eunion 
de facettes. Si de plus $\FC'$ est une facette diff\'erente de $\FC$ et contenue 
dans $\overline{\FC}$, l'assertion (b) montre qu'il existe $\HC \in \AG_2$ 
tel que $\FC' \subseteq \HC$. Donc $\FC' \subseteq \LC_\AG(\FC) \cap \HC$, 
et $\dim \Pos^{-1}\bigl(\LC_\AG(\FC) \cap \HC\bigr) = 
\dim \Pos^{-1}\bigl(\LC_\AG(\FC)\bigr)-1$. D'o\`u (c). 

\medskip

L'assertion (d) d\'ecoule imm\'ediatement de (c). 
\end{proof}

\bigskip

On d\'efinit une relation $\infspe$ (ou $\preccurlyeq_\AG$ s'il est n\'ecessaire de 
pr\'eciser) entre les facettes~: on \'ecrit $\FC \infspe \FC'$ si 
$\overline{\FC} \subseteq \overline{\FC}'$ (c'est-\`a-dire si $\FC \subseteq \overline{\FC}'$). 
La proposition \ref{facettes} (d) montre que~:

\bigskip

\begin{coro}\label{ordre facettes}
La relation $\infspe$ entre les facettes est une relation d'ordre.
\end{coro}

\bigskip

\section{Alg\`ebres de Hecke} 

\medskip

\subsection{Pr\'eliminaires}
Soit $(W,S)$ un groupe de Coxeter ($S$ et, bien s\^ur, $W$ pouvant \^etre 
infinis). Si $s$, $t \in S$, nous \'ecrirons 
$s \sim t$ si $s$ et $t$ sont conjugu\'es dans $W$. Posons $\Sba=S/\!\!\sim$. 
Si $s \in S$, nous noterons $\sba$ sa classe dans $\Sba$. Nous noterons 
$\ZM[\Sba]$ le $\ZM$-module libre de base $\Sba$. 

Soit $\ell : W \to \NM$ la fonction longueur associ\'ee \`a $S$. Si 
$\o \in \Sba$ et si $w \in W$, nous noterons $\ell_\o(w)$ 
le nombre d'apparitions d'\'el\'ements de $\o$ dans une 
expression r\'eduite de $w$ (il est bien connu que cela ne d\'epend 
pas du choix de la d\'ecomposition r\'eduite). Posons 
$$\fonction{\ellb}{W}{\ZM[\Sba]}{w}{\DS{\sum_{\o \in \Sba}} 
\ell_\o(w) \o.}$$

Si $\G$ est un groupe ab\'elien, une application $\ph : W \to \G$ 
est appel\'ee une {\it fonction de poids} si $\ph(ww')=\ph(w)+\ph(w')$ 
pour tous $w$, $w' \in W$ tels que $\ell(ww')=\ell(w)+\ell(w')$. 
Les fonctions $\ell_\o$ (vues comme fonctions \`a valeurs dans $\ZM$) 
et $\ellb$ sont des fonctions de poids. En fait, la fonction $\ellb$ est 
universelle dans le sens suivant~:

\bigskip

\begin{lem}\label{factorisation poids}
Soit $\ph : W \to \G$ une fonction de poids. Alors il existe 
un unique morphisme de groupes $\phba : \ZM[\Sba] \to \G$ 
tel que $\ph = \phba \circ \ellb$. 
\end{lem}

\bigskip

\begin{proof}
Claire.
\end{proof}

\bigskip

Notons $\maps(\Sba,\G)$ l'ensemble des applications $\Sba \to \G$ et 
notons $\weight(W,\G)$ celui des fonctions de poids $W \to \G$. 
Le lemme \ref{factorisation poids} montre qu'il existe des bijections 
canoniques 
$$\weight(W,\G) \stackrel{\sim}{\longleftrightarrow} \maps(\Sba,\G) 
\stackrel{\sim}{\longleftrightarrow} \Hom(\ZM[\Sba],\G).$$
Nous identifierons dans la suite ces trois ensembles. 
Plus pr\'ecis\'ement, nous travaillerons avec des applications 
$\ph : \Sba \to \G$ que nous verrons indiff\'eremment comme 
des morphismes $\ZM[\Sba] \to \G$ ou des fonctions de poids 
$W \to \G$. En particulier, nous pourrons parler aussi bien de 
$\ph(w)$ (pour $w \in W$) que de $\ph(\l)$ (pour $\l \in \ZM[\Sba]$), 
en esp\'erant que cela n'entra\^{\i}ne pas de confusion. 
Par exemple, $\Ker\ph$ est un sous-groupe de $\ZM[\Sba]$ 
(et non pas de $W$~!). 

\bigskip

\subsection{Alg\`ebres de Hecke} 
Fixons tout d'abord les notations en vigueur jusqu'\`a la fin de cet article. 

\medskip

\begin{quotation}
\noindent{\bf Notations.} {\it Soit $\G$ un groupe ab\'elien et 
soit $\ph : \Sba \to \G$ une application.} 
\end{quotation}

\medskip

\noindent Adoptons une notation exponentielle 
pour l'alg\`ebre de groupe de $\G$~: 
$\ZM[\G]=\DS{\mathop{\oplus}_{\g \in \G}} \ZM e^\g$, o\`u $e^\g\cdot e^{\g'}=e^{\g+\g'}$ pour tous $\g$, $\g' \in \G$. Nous noterons 
alors $\HC(W,S,\ph)$ l'{\it alg\`ebre de Hecke de param\`etre $\ph$}, 
c'est-\`a-dire le $\ZM[\G]$-module libre de base $(T_w)_{w \in W}$ muni d'une 
multiplication $\ZM[\G]$-bilin\'eaire 
totalement d\'etermin\'ee par les r\`egles suivantes~:
$$\begin{cases}
T_w T_{w'} = T_{ww'} & \text{si $\ell(ww')=\ell(w)+\ell(w')$,}\\
(T_s - e^{\ph(s)})(T_s+e^{-\ph(s)})=0 & \text{si $s \in S$.}
\end{cases}$$
S'il est n\'ecessaire de pr\'eciser, nous noterons $T_w^\ph$ 
l'\'el\'ement $T_w$ pour rappeler qu'il vit dans l'alg\`ebre de 
Hecke de param\`etre $\ph$. 

Cette alg\`ebre est munie de plusieurs involutions. 
Nous n'utiliserons que la suivante~: si $\g \in \G$ et 
$w\in W$, posons $\overline{e^\g}=e^{-\g}$ et $\overline{T}_w=T_{w^{-1}}^{-1}$ 
(notons que $T_w$ est inversible). Ceci s'\'etend par $\ZM$-lin\'earit\'e 
en une application $\ZM[\G]$-antilin\'eaire $\HC(W,S,\ph) \to \HC(W,S,\ph)$, 
$h \mapsto \hov$ qui est un automorphisme involutif d'anneau. 

\medskip

La construction pr\'ec\'edente est fonctorielle. 
Si $\r : \G \to \G'$ est un morphisme de groupes ab\'eliens, 
alors $\r$ induit une application 
$$\r_* : \HC(W,S,\ph) \longto \HC(W,S,\r \circ \ph)$$
d\'efinie comme \'etant l'unique application $\ZM[\G]$-lin\'eaire 
envoyant $T_w^\ph$ sur $T_w^{\r \circ \ph}$ (ici, $\HC(W,S,\r \circ \ph)$ 
est vu comme une $\ZM[\G]$-alg\`ebre \`a travers le morphisme 
$\ZM[\G] \to \ZM[\G']$ induit par $\r$). Il est alors 
facile de v\'erifier que 
\equat\label{rho morphisme}
\text{\it $\r_*$ est un morphisme de $\ZM[\G]$-alg\`ebres.}
\endequat
Si $h \in \HC(W,S,\ph)$, alors 
\equat\label{rho bar}
\overline{\r_*(h)}=\r_*(\overline{h}).
\endequat
D'autre part, si $\s : \G' \to \G''$ est un autre morphisme 
de groupes ab\'eliens, alors
\equat\label{covariance rho}
(\s \circ \r)_* = \s_* \circ \r_*.
\endequat
Il en d\'ecoule le lemme suivant~:

\bigskip

\begin{lem}\label{rho bij}
Le morphisme de groupes $\r$ est injectif (respectivement surjectif, 
respectivement bijectif) si et 
seulement si le morphisme d'alg\`ebres $\r_*$ l'est.
\end{lem}

\bigskip

\subsection{Alg\`ebre de Hecke g\'en\'erique} 
Soit $R$ l'alg\`ebre de groupe $\ZM[\ZM[\Sba]]$. 
Notons $i : \Sba \to \ZM[\Sba]$ l'application canonique. 
L'alg\`ebre de Hecke $\HC(W,S,i)$ sera alors not\'ee 
${{\mathcal{H}}}(W,S)$~: elle est appel\'ee l'{\it alg\`ebre de 
Hecke g\'en\'erique}. 

Elle est universelle dans le sens que, si on identifie 
l'application $\ph : \Sba \to \G$ avec le morphisme 
$\ph : \ZM[\Sba] \to \G$, l'alg\`ebre de Hecke $\HC(W,S,\ph)$ 
se retrouve \'equip\'ee d'un morphisme canonique de $R$-alg\`ebres 
$\ph_* : \HCB(W,S) \to \HC(W,S,\ph)$, qui fait de $\HC(W,S,\ph)$ une 
sp\'ecialisation de $\HCB(W,S)$. 

\bigskip

\section{Cellules de Kazhdan-Lusztig}

\medskip

\subsection{Base de Kazhdan-Lusztig} 
Pour d\'efinir la base de Kazhdan-Lusztig, nous ferons 
l'hypoth\`ese suivante~: 

\medskip

\begin{quotation}
\noindent{\bf Hypoth\`ese et notations.} {\it Nous supposerons 
jusqu'\`a la fin de cet article que $\G$ est \'equip\'e 
d'un ordre total $\le$ faisant de lui un groupe ab\'elien totalement 
ordonn\'e. Nous noterons 
respectivement $\G_{>0}$, $\G_{\geqslant 0}$, $\G_{< 0}$ et $\G_{\leqslant 0}$ 
l'ensemble des \'el\'ements stritements positifs, positifs ou nuls, 
strictement n\'egatifs et n\'egatifs ou nuls de $\G$.} 
\end{quotation}

\medskip

Si $E$ est un sous-ensemble quelconque de $\G$, nous poserons 
$\ZM[E]=\DS{\mathop{\oplus}_{\g \in E}} \ZM e^\g$. Par exemple, 
$\ZM[\G_{\leqslant 0}]$ est une sous-alg\`ebre de $\ZM[\G]$ dont 
$\ZM[\G_{<0}]$ est un id\'eal. 
Posons 
$$\HC_{<0}(W,S,\ph) = \mathop{\oplus}_{w \in W} \ZM[\G_{< 0}] ~T_w.$$
Alors, si $w \in W$, il existe un unique \'el\'ement 
$C_w \in \HC(W,S,\ph)$ tel que 
$$\begin{cases}
\overline{C}_w = C_w,& \\
C_w \in T_w + \HC_{<0}(W,S,\ph)&
\end{cases}$$
(voir \cite[Theor\`eme 5.2]{lusztig} en g\'en\'eral). Encore une fois, s'il est n\'ecessaire 
de pr\'eciser, l'\'el\'ement $C_w$ sera not\'e $C_w^\ph$. 

\def\pre#1#2{\leqslant_{#1}^{#2}}

La famille $(C_w)_{w \in W}$ forme alors une $\ZM[\G]$-base de 
$\HC(W,S,\ph)$ (appel\'ee {\it base de Kazhdan-Lusztig} \cite[Theorem 5.2]{lusztig}).
Nous noterons $\pre{L}{\ph}$, $\pre{R}{\ph}$ et $\pre{LR}{\ph}$ 
les pr\'eordres d\'efinis par Kazhdan et Lusztig \cite[\S 8.1]{lusztig} 
et nous noterons $\sim_L^\ph$, $\sim_R^\ph$ et $\sim_{LR}^\ph$ les relations 
d'\'equivalence associ\'ees.

Si $w \in W$, et si $? \in \{L, R, LR\}$, 
nous noterons 
$$\Cell_?^\ph(w)=\{x \in W~|~x \sim_?^\ph w\}.$$
Si $?=L$, $R$ ou $LR$, alors $\Cell_?(w)$ 
est appel\'ee la {\it cellule \`a gauche}, la {\it cellule \`a droite} 
ou la {\it cellule bilat\`ere} de $w$ (pour $(W,S,\ph)$). 

\bigskip

\subsection{Repr\'esentations cellulaires} 
Si $C$ est une cellule \`a gauche pour $(W,S,\ph)$, on peut lui associer un 
$\HC(W,S,\ph)$-module $M_C^\ph$ \cite[\S 8.3]{lusztig}. Notons $\aug : \ZM[\G] \to \ZM$ le morphisme 
d'augmentation. Voyons $\ZM$ (ou $\QM$) comme une $\ZM[\G]$-alg\`ebre 
\`a travers $\aug$, on a un isomorphisme d'anneaux 
$\ZM \otimes_{\ZM[\G]} \HC(W,S,\ph) \simeq \ZM W$. Notons 
$\ZM M_C^\ph$ (resp. $\QM M_C^\ph$) le $\ZM W$-module (resp. $\QM W$-module) 
$\ZM \otimes_{\ZM[\G]} M_C^\ph$ (resp. $\QM \otimes_{\ZM[\G]} M_C^\ph$). 

Si $W$ est {\it fini}, nous noterons 
$\chi_C^\ph$ le caract\`ere de $\QM M_C^\ph$. 

\bigskip

\subsection{Morphismes strictement croissants} 
L'objectif principal de cet article est d'\'etudier 
comment se comporte la partition en cellules lorsque le triplet 
$(\G,\le,\ph)$ varie. La premi\`ere remarque facile 
est que cette partition ne change pas si on compose 
$\ph$ avec un morphisme strictement croissant:

\bigskip

\begin{prop}\label{strictement croissant}
Soit $\G'$ un groupe ab\'elien totalement ordonn\'e, 
soit $\r : \G \to \G'$ un morphisme de groupes strictement croissant et 
soit $? \in \{L,R,LR\}$. 
Alors~:
\begin{itemize}
 \itemth{a} $\r_*$ est injectif.

\itemth{b} Si $w \in W$, $\r_*(C_w^\ph)=C_w^{\r \circ \ph}$. 

\itemth{c} Les relations 
$\le_?^\ph$ et $\le_?^{\r \circ \ph}$ sont \'egales (de m\^eme 
que les relations $\sim_?^\ph$ et $\sim_?^{\r \circ \ph}$). 

\itemth{d} Si $w \in W$, alors  $\Cell_?^\ph(w)=\Cell_?^{\r \circ \ph}(w)$. 
\end{itemize}
\end{prop}

\begin{proof}
L'injectivit\'e de $\r_*$ d\'ecoule du fait que $\r$ lui-m\^eme 
est injectif (voir le lemme \ref{rho bij}). 
D'o\`u (a). Montrons (b). Posons $C=\r_*(C_w^\ph)$. 
Alors, d'apr\`es \ref{rho bar}, on a $\Cov=C$. D'autre part, 
la stricte croissance de $\r$ implique que 
\equat\label{rho positif}
\r_*(\HC_{<0}(W,S,\ph)) \subseteq \HC_{< 0}(W,S,\r \circ \ph).
\endequat
Par cons\'equent, 
$C - T_w^{\r \circ \ph} \in \HC_{<0}(W,S,\r \circ \ph)$. 
Donc $C=C_w^{\r \circ \ph}$ d'apr\`es la caract\'erisation des
\'el\'ements de la base de Kazhdan-Lusztig. 

Montrons maintenant (c). Si $x$, $y$ et $z \in W$, alors 
il r\'esulte de (b) que 
\equat\label{rho h}
\r(h_{x,y,z}^\ph)=h_{x,y,z}^{\r \circ \ph}
\endequat
(il suffit d'appliquer $\r_*$ \`a la relation d\'efinissant 
les \'el\'ements $h_{x,y,z}$). L'assertion (c) d\'ecoule 
alors de cette observation, du fait que 
$\r : \ZM[\G] \to \ZM[\G']$ est injectif et de la d\'efinition 
des relations $\le_?^\ph$ et $\le_?^{\r \circ \ph}$.
 
L'assertion (d) d\'ecoule imm\'ediatement de (c).
\end{proof}

\bigskip

\subsection{Parties positives de ${\boldsymbol{\ZM[\Sba]}}$} 
Soit $X$ une partie positive de $\ZM[\Sba]$. Posons 
$\G_X = \ZM[\Sba]/(X \cap (-X))$ et notons $\le_X$ l'ordre total 
sur $\G_X$ d\'efini dans \ref{ordre X}. Notons $\ph_X : \Sba \to \G_X$ 
l'application canonique. Pour simplifier, les relations 
$\le_?^{\ph_X}$ et $\sim_?^{\ph_X}$ seront not\'ees 
$\le_?^X$ et $\sim_?^X$. De m\^eme, si $w \in W$, 
nous noterons $\Cell_?^X(w)$ l'ensemble $\Cell_?^{\ph_X}(w)$ et, 
si $C$ est une cellule pour $(W,S,\ph_X)$, nous noterons 
$\ZM M_C^X$ le $\ZM W$-module $\ZM M_C^{\ph_X}$ 
(et, si $W$ est {\it fini}, nous noterons $\chi_C^X$ le caract\`ere 
$\chi_C^{\ph_X}$).

La proposition suivante montre que la famille de triplets 
$\bigl((\G_X,\le_X,\ph_X)\bigr)_{X \in \positif(\ZM[\Sba])}$ 
est essentiellement exhaustive~:

\bigskip

\begin{prop}
Posons $X=\Pos(\ph)$. Alors il existe un unique morphisme de 
groupes $\phba : \G_X \to \G$ tel que $\ph = \phba \circ \ph_X$. 
Ce morphisme $\phba$ est strictement croissant. 
\end{prop}

\begin{proof} 
En effet, $\Ker \ph = X \cap (-X)$ donc l'existence et l'unicit\'e de 
$\phba$ est assur\'ee. L'injectivit\'e de $\phba$ est imm\'ediate. 
D'autre part, si $\g$, $\g' \in \G_X$ sont 
tels que $\g \le \g'$ et si $\l \in \ZM[\Sba]$ est un repr\'esentant de 
$\g'-\g$, alors $\l \in X$ d'apr\`es \ref{relation independante}. 
Donc $\ph(\l) \ge 0$. En d'autres termes $\phba(\g'-\g) \ge 0$, 
c'est-\`a-dire $\phba(\g) \le \phba(\g')$. Donc $\phba$ est 
croissant~: la croissance stricte r\'esulte de l'injectivit\'e.
\end{proof}

\bigskip

\begin{coro}\label{phi X}
Posons $X=\Pos(\ph)$ et soit $? \in \{L,R,LR\}$. Alors~:
\begin{itemize}
\itemth{a} Les relations 
$\le_?^\ph$ et $\le_?^X$ sont \'egales (de m\^eme 
que les relations $\sim_?^\ph$ et $\sim_?^X$). 

\itemth{b} Si $w \in W$, alors 
$\Cell_?^\ph(w)=\Cell_?^X(w)$.

\itemth{c} Si $C=\Cell_L^\ph(w)$, alors $\ZM M_C^\ph \simeq \ZM M_C^X$. 
\end{itemize}
\end{coro}

\bigskip

\remark{phi X traduction}
Soient $\l \in \ZM[\Sba]$ et $X \in \positif(\ZM[\Sba])$. On a alors les \'equivalences suivantes~:
\begin{itemize}
\itemth{a} $X \in \UC(\l)$ si et seulement si $\ph_X(\l) < 0$.
\itemth{b} $X \in \HC_\l$ si et seulement si $\ph_X(\l) = 0$.
\itemth{c} $X \in \overline{\UC(\l)}$ si et seulement si $\ph_X(\l) \le 0$.\finl
\end{itemize}

\bigskip

\subsection{Changement de signe}
Avant d'\'enoncer quelques conjectures sur le comportement 
des relations $\sim_?^X$ lorsque $X$ varie dans $\positif(\ZM[\Sba])$, 
nous allons \'etudier l'effet du changement de signes de certaines 
valeurs de $\ph$. Nous travaillerons ici sous l'hypoth\`ese suivante~:

\medskip

\begin{quotation}
\noindent{\it 
Soit $S=S_+ \hskip1mm\dot{\cup}\hskip1mm S_-$ une partition de $S$ telle que, si $s \in S_+$ et $t \in S_-$, 
alors $s$ et $t$ ne sont pas conjugu\'es dans $W$. Soit $\ph' : \Sba \to \G$ 
l'application d\'efinie par
$$\ph'(\sba)=\begin{cases}
\ph(\sba) & \text{si $s \in S_+$,}\\
-\ph(\sba) & \text{si $s \in S_-$.}
\end{cases}$$}
\end{quotation}

\medskip

Si $w \in W$, on pose 
$$\ell_\pm(w)=\sum_{\o \in \Sba_\pm} \ell_\o(w),$$
de sorte que $\ell(w)=\ell_+(w) + \ell_-(w)$. On note 
$\th : \HC(W,S,\ph) \longto \HC(W,S,\ph')$ l'unique application 
$\ZM[\G]$-lin\'eaire telle que
$$\th(T_w^\ph)=(-1)^{\ell_-(w)} T_w^{\ph'}.$$
Un calcul \'el\'ementaire montre que~:

\bigskip

\begin{prop}\label{signe}
L'application $\th : \HC(W,S,\ph) \longto \HC(W,S,\ph')$ 
est un isomorphisme de $\ZM[\G]$-alg\`ebres. De plus, si $h \in \HC(W,S,\ph)$, 
alors $\th(\overline{h})=\overline{\th(h)}$. Par cons\'equent, si 
$w \in W$, alors
$$\th(C_w^\ph) = (-1)^{\ell_-(w)} C_w^{\ph'}.$$
\end{prop}

\bigskip

\begin{coro}\label{signe cellules}
Si $? \in \{L,R,LR\}$, alors les relations $\pre{?}{\ph}$ et 
$\pre{?}{\ph'}$ co\"\i ncident. De m\^eme, les relations 
$\sim_?^\ph$ et $\sim_?^{\ph'}$ co\"\i ncident.
\end{coro}

\bigskip

Si $\o \in \Sba$, on note $\t_\o$ la sym\'etrie $\ZM$-lin\'eaire 
sur $\ZM[\Sba]$ telle que $\t_\o(\o)=-\o$ et $\t_\o(\o')=\o'$ si 
$\o' \neq \o$. C'est un automorphisme de $\ZM[\Sba]$~: il induit 
donc un hom\'eomorphisme $\t_\o^*$ de $\positif(\ZM[\Sba])$. 

\bigskip

\begin{coro}\label{symetrie}
Si $\o \in \Sba$, si $X \in \positif(\ZM[\Sba])$ et si 
$? \in \{L,R,LR\}$, alors les relations $\sim_?^X$ et 
$\sim_?^{\t_\o^*(X)}$ co\"{\i}ncident.
\end{coro}

\bigskip

\begin{proof}
L'application $\t_\o$ induit un isomorphisme strictement 
croissant 
$$\t_\o : \G_{\t_\o^*(X)} \stackrel{\sim}{\longrightarrow} \G_X.$$
Par cons\'equent, les relations $\sim_?^{\ph_{\t_\o^*(X)}}$ 
et $\sim_?^{\t_o \circ \ph_{\t_\o^*(X)}}$ co\"{\i}ncident 
(voir la proposition \ref{strictement croissant}). 
Posons $\ph_X'=\t_o \circ \ph_{\t_\o^*(X)} : \Sba \to \G_X$. 
Or,
$$\ph_X'(s)=\begin{cases}
\ph_X(s) & \text{si $s \not\in \o$,}\\
-\ph_X(s) & \text{si $s \in \o$.}
\end{cases}$$
Donc le r\'esultat d\'ecoule du corollaire \ref{signe cellules}.
\end{proof}

\bigskip

Notons $\g : W \to \{1,-1\}$ l'unique caract\`ere lin\'eaire de $W$ tel que 
$$\g(s)=
\begin{cases}
1 & \text{si $s \in S_+$},\\
-1 & \text{si $s \in S_-$}.
\end{cases}$$
Notons $\ZM_\g$ le $\ZM W$-module (irr\'eductible) de dimension $1$ sur lequel 
$W$ agit via le caract\`ere $\g$. Si $C$ est une cellule \`a gauche pour $(W,S,\ph)$, 
alors un calcul \'el\'ementaire utilisant l'isomorphisme $\th$ montre que
\equat\label{module signe}
\ZM M_C^\ph \simeq \ZM_\g \otimes_\ZM \ZM M_C^{\ph'}.
\endequat

\bigskip

\subsection{Param\`etres nuls}
Il r\'esulte du corollaire \ref{signe cellules} que le calcul des cellules de 
Kazhdan-Lusztig peut se ramener au cas o\`u $\ph$ est \`a valeurs dans 
$\G_{\ge 0}$. Nous allons \'etudier ici ce qu'il se passe lorsque 
certains param\`etres sont nuls. 
Nous travaillerons sous les hypoth\`eses suivantes~:

\medskip

\begin{quotation}
\noindent{\it Dans cette sous-section, et dans cette sous-section seulement, nous fixons 
une partition $S=I \hskip1mm\dot{\cup}\hskip1mm J$ de $S$ telle que, si $s \in I$ et $t \in J$, 
alors $s$ et $t$ ne sont pas conjugu\'es dans $W$. On note $W_I$ 
le sous-groupe de $W$ engendr\'e par $I$ et on pose 
$$\Jti=\{wtw^{-1}~|~w \in W_I,~t \in J\}.$$
Soit $\Wti$ le sous-groupe de $W$ engendr\'e par $\Jti$. 
Nous supposerons de plus que, 
$$\text{\it si $s \in I$, alors $\ph(s)=0$.}\leqno{\hskip1.3cm(*)}$$}
\end{quotation}

D'apr\`es \cite[th\'eor\`eme 1]{bonnafe dyer}, $(\Wti,\Jti)$ est 
un groupe de Coxeter et
\equat\label{dyer bis}
W=W_I \ltimes \Wti
\endequat
Si $\tti \in \Jti$, on notera $\nu(\tti)$ l'unique \'el\'ement de 
$J$ tel que $\tti$ soit conjugu\'e \`a $\nu(\tti)$ 
(voir \cite{bonnafe dyer}). On pose
$$\pht(\tti)=\ph(\nu(\tti)).$$
Il r\'esulte de \cite[(4)]{bonnafe dyer} que, si 
$\tti$ et $\tti'$ sont deux \'el\'ements de $\Jti$ qui sont conjugu\'es 
sous $\Wti$, alors 
\equat\label{pht}
\pht(\tti)=\pht(\tti').
\endequat
Cela montre que l'on peut d\'efinir une $\ZM[\G]$-alg\`ebre de Hecke 
$\HC(\Wti,\Jti,\pht)$. Le groupe $W_I$ agit sur $\Wti$ en stabilisant 
$\Jti$ et la fonction $\pht$, donc il agit naturellement sur 
l'alg\`ebre de Hecke $\HC(\Wti,\Jti,\pht)$. On peut alors former 
le produit semi-direct
$$W_I \ltimes \HC(\Wti,\Jti,\pht).$$
C'est une $\ZM[\G]$-alg\`ebre de $\ZM[\G]$-base $(x \cdot T_w^\pht)_{x \in W_I,w \in \Wti}$. 
Pour finir, posons
$$\HCt=\mathop{\oplus}_{w \in \Wti} \ZM[\G]~T_w^\ph \quad \subseteq \HC(W,S,\ph).$$

\begin{prop}\label{semi-direct}
Rappelons que $\ph(s)=0$ si $s \in I$. Alors~:
\begin{itemize}
\itemth{a} $\HCt$ est une sous-alg\`ebre de $\HC(W,S,\ph)$.

\itemth{b} L'unique application $\ZM[\G]$-lin\'eaire 
$\th : W_I \ltimes \HC(\Wti,\Jti,\pht) \longto \HC(W,S,\ph)$ qui envoie 
$x \cdot T_w^\pht$ sur $T_{xw}^\ph$ ($x \in W_I$, $w \in \Wti$) 
est un isomorphisme de $\ZM[\G]$-alg\`ebres. Il envoie $\HC(\Wti,\Jti,\pht)$ 
isomorphiquement sur $\HCt$.

\itemth{c} Si $x \in W_I$ et $w \in W$, alors 
$C_x^\ph=T_x^\ph$, $C_{xw}^\ph=T_x^\ph C_w^\ph$, 
$C_{wx}^\ph=C_w^\ph T_x^\ph$. Si de plus 
$w \in \Wti$, alors $\th(C_w^\pht)=C_w^\ph$.
\end{itemize}
\end{prop}

\begin{proof}
Puisque $\ph(I)=\{0\}$, on a, pour tout $s \in I$, 
$$(T_s^\ph)^2 = 1\quad\text{et}\quad \overline{T_s^\ph} = T_s^\ph.$$
Cela montre que, si $x \in W_I$ et $w \in W$, alors 
$T_x^\ph T_w^\ph=T_{xw}^\ph$ et $T_w^\ph T_x^\ph=T_{wx}^\ph$. 
On en d\'eduit que $C_x^\ph=T_x^\ph$, $C_{xw}^\ph=T_x^\ph C_w^\ph$, 
$C_{wx}^\ph=C_w^\ph T_x^\ph$. C'est la premi\`ere assertion du (c). 

\medskip

D'autre part, si $x \in W_I$ et $t \in J$, alors 
$$T_{xtx^{-1}}^\ph=T_x^\ph T_t^\ph (T_x^\ph)^{-1}$$
et donc, en posant $\tti=xtx^{-1}$, on a 
$$(T_\tti^\ph - e^{\pht(\tti)})(T_{\tti}^\ph+e^{-\pht(\tti)})=0.$$
Pour montrer les assertions (a) et (b), il ne reste qu'\`a montrer 
l'assertion suivante~: si $w$ et $w'$ sont deux \'el\'ements de 
$\Wti$ tels que $\ellt(ww')=\ellt(w)+\ellt(w')$ (ici, $\ellt$ d\'esigne 
la fonction longueur sur $\Wti$ associ\'ee \`a $\Jti$), alors 
$$T_w^\ph T_{w'}^\ph =T_{ww'}^\ph.\leqno{(?)}$$
Un raisonnement par r\'ecurrence (sur $\ellt(w)$) 
\'el\'ementaire permet de se ramener au cas o\`u $\ellt(w)=1$, c'est-\`a-dire 
$w \in \Jti$. Dans ce cas, $w=xtx^{-1}$, avec $t \in J$ et $x \in W_I$. 
La fonction $\ellt$ \'etant $W_I$-invariante, on a 
$\ellt(t x^{-1}w'x)=1+\ellt(x^{-1}w' x)$. Donc, 
d'apr\`es le \cite[corollaire 7]{bonnafe dyer}, 
on a $\ell(t x^{-1}w'x)=1+\ell(x^{-1}w'x)$. Il en r\'esulte que 
$$T_t^\ph T_{x^{-1}w'x}^\ph=T_{tx^{-1}w'x}^\ph.$$
On obtient alors (?) en conjuguant cette \'egalit\'e par $T_x^\ph$. 
D'o\`u (a) et (b). 

\medskip

Pour montrer la derni\`ere assertion de (c), il suffit de remarquer que, 
si $h \in \HC(\Wti,\Jti,\pht)$, alors $\th(\overline{h})=\overline{\th(h)}$, 
et d'utiliser la caract\'erisation de la base de Kazhdan-Lusztig. 
\end{proof}

On en d\'eduit imm\'ediatement le corollaire suivant~:

\begin{coro}\label{KL semi-direct}
Supposons que $\ph(I) =\{0\}$. 
Alors les cellules \`a gauche (resp. \`a droite, resp. bilat\`eres) pour $(W,S,\ph)$
sont de la forme $W_I \cdot C$ (resp. $C \cdot W_I$, resp. $W_I \cdot C \cdot W_I$), 
o\`u $C$ est une cellule \`a gauche (resp. droite, resp. bilat\`ere) 
pour $(\Wti,\Jti,\pht)$.
\end{coro}

\bigskip

\section{Conjectures}

\medskip
 
Nos conjectures portent sur le comportement des relations $\sim_?^X$ lorsque 
$X$ varie dans l'espace topologique $\positif(\ZM[\Sba])$. 
On notera $\RM[\Sba]$ le $\RM$-espace vectoriel $\RM \otimes_\ZM \ZM[\Sba]$. 
Un arrangement d'hyperplans rationnels 
$\AG$ de $\positif(\ZM[\Sba])$ sera dit {\it complet} si $\HC_\o \in \AG$ pour tout 
$\o \in \Sba$ (rappelons que $\HC_\o$ est l'ensemble des 
$X \in \positif(\ZM[\Sba])$ tels que $\o \in X \cap -X$). 
Si $\AG$ est un arrangement complet et si $\FC$ est une $\AG$-facette, 
nous noterons $W_\FC$ le sous-groupe parabolique standard 
de $W$ engendr\'e par la r\'eunion des 
orbites $\o \in \Sba$ telles que $\FC \subseteq \HC_\o$.

\bigskip

\noindent{\bf Conjecture A.} 
{\it Il existe un arrangement complet fini $\AG$ d'hyperplans rationnels 
de $\positif(\ZM[\Sba])$ v\'erifiant les propri\'et\'es suivantes (pour $? \in \{L,R,LR\}$)~:
\begin{itemize}
 \itemth{a} Si $X$ et $Y$ sont deux parties positives de $\ZM[\Sba]$ 
appartenant \`a la m\^eme $\AG$-facette $\FC$, alors les relations $\sim_?^X$ et 
$\sim_?^Y$ co\"{\i}ncident. Nous la noterons 
$\sim_?^\FC$.

\itemth{b} Soit $\FC$ une $\AG$-facette. Alors les cellules (i.e. classes d'\'equivalence) 
pour la relation $\sim_?^\FC$ sont les parties minimales $C$ de $W$ satisfaisant 
aux conditions suivantes~:
\begin{itemize}
\itemth{b1} Pour toute chambre $\CC$ telle que $\FC \infspe \CC$, $C$ est une r\'eunion 
de cellules pour la relation $\sim_?^\CC$~;

\itemth{b2} $C$ est stable par translation par $W_\FC$ (\`a gauche si $?=L$, \`a droite 
si $?=R$, \`a gauche et \`a droite si $?=LR$).
\end{itemize}
\end{itemize}}

\bigskip

Si la conjecture $A$ est vraie pour $(W,S)$ et deux arrangements 
complets finis $\AG$ et $\AG'$ d'hyperplans rationnels de 
$\positif(\ZM[\Sba])$, alors elle est vraie pour $(W,S)$ 
et l'arrangement complet fini $\AG \cap \AG'$. Cela montre que, si 
la conjecture A est vraie pour $(W,S)$, il existe un arrangement 
complet fini 
d'hyperplans rationnels qui est minimal et pour lequel les assertions 
de la conjecture A sont vraies : nous les appellerons 
{\it hyperplans essentiels} (en r\'ef\'erence aux travaux de M. Chlouveraki 
sur les blocs de Rouquier des alg\`ebres de Hecke 
cyclotomiques \cite[\S 4.4]{chlouveraki}) de $(W,S)$). 

Le corollaire \ref{symetrie} montre que, si la conjecture 
A est vraie pour $(W,S)$, alors l'ensemble des hyperplans 
essentiels de $(W,S)$ est stable par l'action de toutes les 
sym\'etries $\t_\o^*$ ($\o \in \Sba$). 

\bigskip

\noindent{\sc Remarque - } 
Fixons $w \in W$. La condition (a) de la conjecture A dit que l'application 
$\positif(\L) \to \PC(W)$, $X \mapsto \Cell_?^X(w)$ est constante 
sur les facettes. La condition (b1) dit que cette m\^eme application est 
semi-continue sup\'erieurement (rappelons qu'une application $f : \BC \to B$, 
o\`u $\BC$ est un espace topologique et $B$ est un ensemble ordonn\'e, 
est dite {\it semi-continue sup\'erieurement} si, pour tout $b \in B$, 
l'ensemble $\{x \in \BC~|~f(x) < b\}$ est ouvert). 
Notons pour finir que toute cellule pour $\sim_\FC^?$ est forc\'ement stable par translation 
par $W_\FC$ (voir le corollaire \ref{KL semi-direct}~: en effet, 
si $X \in \FC$ et $\FC \subseteq \HC_\o$, alors $\ph_X(s)=0$ pour 
tout $s \in \o$)~: cela justifie la propri\'et\'e (b2). 

En revanche, la condition de minimalit\'e dans (b) 
apporte une pr\'ecision \'etonnante~: pour conna\^{\i}tre les cellules de 
$W$ pour un choix de fonction $\ph$, il suffirait de conna\^{\i}tre 
les relations $\sim_?^\CC$, lorsque $\CC$ est une chambre.\finl

\bigskip

%
%
Il semble raisonnable d'imaginer que la Conjecture A est compatible avec 
la construction des repr\'esentations cellulaires~:

\bigskip

\noindent{\bf Conjecture B.} 
{\it Supposons la conjecture A vraie pour $(W,S)$ et notons 
$\AG$ l'ensemble des hyperplans essentiels de $(W,S)$. Soit 
$X \in \positif(\ZM[\Sba])$, soit $C$ une cellule \`a gauche 
pour $(W,S,\ph_X)$ et soit $\CC$ une $\AG$-chambre dans 
$\positif(\ZM[\Sba])$ telle que $X \in \overline{\CC}$. Soit 
$Y \in \CC$. D'apr\`es la conjecture A, il existe des cellules 
\`a gauche $C_1$,\dots, $C_d$ pour $(W,S,\ph_Y)$ telles que 
$C = C_1 \hskip1mm\dot{\cup}\hskip1mm C_2 \hskip1mm\dot{\cup}\hskip1mm 
\cdots \hskip1mm\dot{\cup}\hskip1mm C_d$. Alors il existe une filtration 
$M_0=0 \subseteq M_1 \subseteq \cdots \subseteq M_d=\ZM M_C^X$ 
du $\ZM W$-module $\ZM M_C^X$ et une permutation $\s \in \SG_d$ 
telle que 
$$M_i/M_{i-1} \simeq \ZM M_{C_{\s(i)}}^Y$$}

\bigskip

Une version (beaucoup) plus faible est donn\'ee par~:

\bigskip

\noindent{\bf Conjecture $\Bb^{\boldsymbol{-}}$.} 
{\it Supposons $W$ fini. 
Supposons la conjecture A vraie pour $(W,S)$ et notons $\AG$ l'ensemble 
des hyperplans essentiels de $(W,S)$. Soit $X \in \positif(\ZM[\Sba])$, 
soit $C$ une cellule \`a gauche pour $(W,S,\ph_X)$ et soit $\CC$ 
une $\AG$-chambre dans $\positif(\ZM[\Sba])$ telle que 
$X \in \overline{\CC}$. Soit $Y \in \CC$. D'apr\`es la conjecture A, 
il existe des cellules \`a gauche $C_1$,\dots, $C_d$ 
pour $(W,S,\ph_Y)$ telles que $C = C_1 \hskip1mm\dot{\cup}\hskip1mm C_2 
\hskip1mm\dot{\cup}\hskip1mm \cdots \hskip1mm\dot{\cup}\hskip1mm C_d$. 
Alors 
$$\chi_C^X=\sum_{i=1}^d \chi_{C_i}^Y.$$}

\bigskip

\section{Exemples}

\bigskip

\subsection{Groupes di\'edraux finis}
Supposons dans cette sous-section, et dans cette sous-section seulement, 
que $S=\{s,t\}$ et que $st$ est d'ordre {\it fini} $2m$ avec $m \ge 2$. 
Notons $w_0=(st)^m=(ts)^m$ l'\'el\'ement le plus long de $W$~: il est central.
Si $w \in W$, on pose $\RC(w)=\{u \in S~|~wu < w\}$. Posons
$$C_s=\{w \in W~|~\RC(w)=\{s\}\}\quad\text{et}\quad
C_t=\{w \in W~|~\RC(w)=\{t\}\}.$$
Un calcul facile \cite[\S 8.7]{lusztig} montre que la partition de $W$ en 
cellules \`a gauche pour $(W,S,\ph)$ est donn\'ee par le tableau suivant 
(lorsque $\ph$ est \`a valeurs dans $\G_{\geqslant 0}$)~:

\bigskip

$$\begin{array}{|c|c|}
\hline
\espace \ph & \text{Cellules \`a gauche} \\
\hline
\hline
\espace 0=\ph(s)=\ph(t) & W \\
\espace 0=\ph(s) < \ph(t) & \{1,s\},~C_s\setminus\{s\},~ 
C_t \setminus\{sw_0\},~ \{sw_0,w_0\} \\
\espace 0 < \ph(s) < \ph(t) & \{1\},~ \{s\},~ C_s \setminus\{s\},~ 
C_t \setminus\{sw_0\},~
\{sw_0\},~ \{w_0\} \\
\espace 0 < \ph(s)=\ph(t) & \{1\},~ C_s,~C_t,~\{w_0\} \\
\espace 0 < \ph(t) < \ph(s) & \{1\},~ \{t\},~ 
C_s\setminus\{tw_0\},~ C_t\setminus\{t\},~ \{tw_0\},~ 
\{w_0\} \\
\espace 0=\ph(t) < \ph(s) & \{1,t\},~ 
C_s\setminus\{tw_0\},~ C_t\setminus\{t\},~ \{tw_0,w_0\}\\
\hline
\end{array}$$

\bigskip

Par cons\'equent~:

\bigskip

\begin{prop}\label{diedral}
Les conjectures A et $B^-$ sont vraies lorsque $|S|=2$ et $|W| < \infty$. 
Les hyperplans essentiels de $(W,S)$ sont $\HC_s$, $\HC_t$, $\HC_{s-t}$ et 
$\HC_{s+t}$. 
\end{prop}

\begin{proof}
L'\'enonc\'e des conjectures faisant appara\^{\i}tre la topologie 
exotique sur $\positif(\ZM[\Sba])$, nous allons ici d\'etailler le passage 
de la table pr\'ec\'edente \`a la preuve des conjectures A et B. 

Tout d'abord, on identifie $\Sba$ \`a $S$, ce qui permet d\'ecrire
$$\ZM[\Sba] = \ZM s \oplus \ZM t.$$
Si on note $\AG=\{\HC_s,\HC_t,\HC_{s-t},\HC_{s+t}\}$, alors $\AG$ est un 
arrangement complet fini d'hyperplans rationnels de $\positif(\ZM[\Sba])$. 
Dessinons l'image inverse dans $\RM[\Sba]^*$ (sous 
l'application $\Pos$) de l'arrangement d'hyperplans et 
de ses facettes. Pour cela, notons $(s^*,t^*)$ la base 
duale de la base $(s,t)$ de $\RM[\Sba]$. 
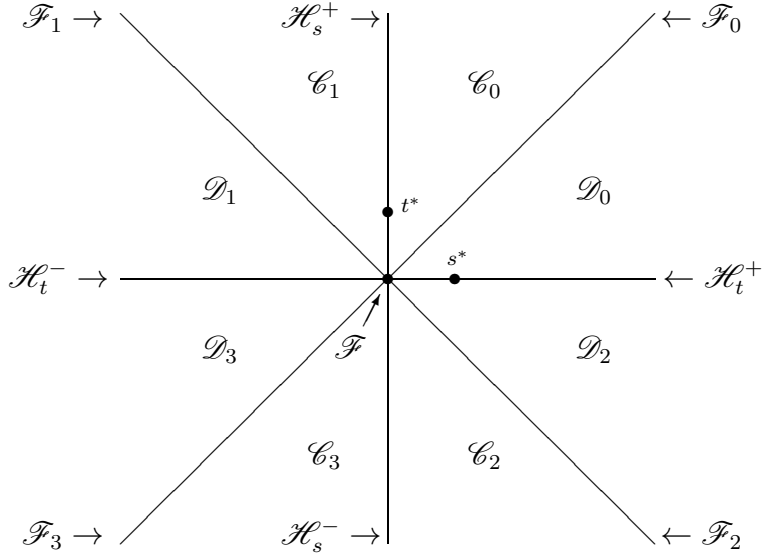
\medskip

\begin{center}
\begin{picture}(200,200)
\put(0,100){\line(1,0){200}}
\put(0,0){\line(1,1){200}}
\put(100,0){\line(0,1){200}}
\put(0,200){\line(1,-1){200}}
\put(130,170){$\CC_0$}
\put(170,130){$\DC_0$}
\put(70,170){$\CC_1$}
\put(30,130){$\DC_1$}
\put(130,30){$\CC_2$}
\put(170,70){$\DC_2$}
\put(70,30){$\CC_3$}
\put(30,70){$\DC_3$}
\put(100,100){\circle*{4}}
\put(125,100){\circle*{4}}\put(122,105){$\SS{s^*}$}
\put(100,125){\circle*{4}}\put(105,125){$\SS{t^*}$}
\put(-42,97){$\HC_t^- \to$}
\put(-37,195){$\FC_1 \to$}
\put(-37,0){$\FC_3 \to$}
\put(202,195){$\leftarrow\FC_0$}
\put(60,195){$\HC_s^+ \to$}
\put(60,0){$\HC_s^- \to$}
\put(202,0){$\leftarrow\FC_2$}
\put(203,97){$\leftarrow\HC_t^+$}
\put(91,82){\vector(1,2){6}}
\put(80,70){$\FC$}
\end{picture}
\end{center}

Les $\AG$-facettes sont donc
$$\FC=\{\ZM[\Sba]\},\quad \HC_s^\pm,\quad \HC_t^\pm, \quad \FC_i,\quad\CC_i,\quad
\DC_i\quad (0 \le i \le 3).$$
Si l'on veut une d\'efinition formelle, voici quelques exemples, 
qui montrent que l'on a bien d\'efini ainsi des facettes~:
$$\FC=\HC_s \cap \HC_t \cap \HC_{s-t} \cap \HC_{s+t},$$
$$\HC_t^+=\HC_t \cap \UC(-s) = \HC_t \cap \UC(-s) \cap \UC(t-s) \cap \UC(-t-s),$$
$$\FC_0=\UC(-s) \cap \HC_{s-t} = \UC(-s) \cap \UC(-t) \cap \HC_{s-t} \cap \UC(-s-t),$$
$$\CC_0=\UC(-s) \cap \UC(s-t)=\UC(-s) \cap \UC(-t) \cap \UC(s-t) \cap \UC(-s-t).
\leqno{\text{et}}$$
Soit $X \in \positif(\ZM[\Sba])$. D'apr\`es la remarque \ref{phi X traduction}, on a les \'equivalences suivantes~:
\begin{itemize}
\itemth{a} $0=\ph_X(s)=\ph_X(t)$ si et seulement si $X \in \FC$.

\itemth{b} $0 = \ph_X(s) < \ph_X(t)$ si et seulement si $X \in \HC_s^+$. 

\itemth{c} $0 < \ph_X(s) < \ph_X(t)$ si et seulement si $X \in \CC_0$.

\itemth{d} $0 < \ph_X(s) = \ph_X(t)$ si et seulement si $X \in \FC_0$.

\itemth{e} $0 < \ph_X(t) < \ph_X(s)$ si et seulement si $X \in \DC_0$.

\itemth{f} $0 =\ph_X(t) < \ph_X(s)$ si et seulement si $X \in \HC_t^+$.
\end{itemize}
En utilisant la table pr\'ec\'edente, les assertions (a)--(f) ci-dessus et les sym\'etries 
$\t_s^*$ et $\t_t^*$ (voir Corollaire \ref{symetrie}), on obtient la table 
suivante qui donne la r\'epartitions en cellules \`a gauche pour $(W,S,X)$ 
selon la facette \`a laquelle appartient $X$~:
\bigskip

$$\begin{array}{|c|c|}
\hline
\espace X \in ? & \text{Cellules \`a gauche} \\
\hline
\hline
\espace \{\ZM[\Sba]\} & W \\
\espace \HC_s^\pm & \{1,s\},~C_s\setminus\{s\},~ 
C_t \setminus\{sw_0\},~ \{sw_0,w_0\} \\
\espace \CC_i~(1 \le i \le 4) & \{1\},~ \{s\},~ 
C_s \setminus\{s\},~ C_t \setminus\{sw_0\},~
\{sw_0\},~ \{w_0\} \\
\espace \FC_i~(1 \le i \le 4) & \{1\},~ C_s,~C_t,~\{w_0\} \\
\espace  \DC_i~(1 \le i \le 4) & \{1\},~ \{tw_0\},~ C_s\setminus\{tw_0\},~ 
C_t\setminus\{t\},~ \{t\},~ \{w_0\} \\
\espace \HC_t^\pm & \{1,t\},~ C_s\setminus\{tw_0\},~ 
C_t\setminus\{t\},~ \{tw_0,w_0\}\\
\hline
\end{array}$$

\bigskip

On peut alors v\'erifier sur cette table les faits suivants~:
\begin{itemize}
\item Si $X$ et $X'$ appartiennent \`a la m\^eme facette, alors les relations 
$\sim_L^X$ et $\sim_L^{X'}$ co\"{\i}ncident.

\item Si $X\in \FC_i$, alors les seules chambres $\CC$ telles que 
$X \in \overline{\CC}$ sont $\CC_i$ et $\DC_i$~: il est alors facile de 
voir que $\sim_L^{\FC_i}$ est effectivement engendr\'ee par $\sim_L^{\CC_i}$ et 
$\sim_L^{\DC_i}$.

\item Si $X \in \HC_s^+$, alors les seules chambres $\CC$ telles que $X \in \overline{\CC}$ 
sont $\CC_0$ et $\CC_1$~: il est alors facile de voir qu'une cellule \`a gauche pour $(W,S,\ph_X)$ est bien 
une partie minimale de $W$ qui est stable par translation \`a gauche par $<s>$ 
tout en \'etant r\'eunion de cellules \`a gauche pour $(W,S,\CC_0)$ (et r\'eunion de 
cellules \`a gauche pour $(W,S,\CC_1)$, mais, 
par sym\'etrie, ce sont aussi des cellules \`a gauche pour $(W,S,\CC_0)$). 
\end{itemize}
Pour obtenir les autres \'enonc\'es de la conjecture A concernant les cellules 
\`a gauche, on peut alors utiliser les sym\'etries (via le corollaire \ref{symetrie}) 
et l'automorphisme de $W$ \'echangeant $s$ et $t$. Les \'enonc\'es correspondants 
pour les cellules \`a droite ou bilat\`eres se montrent de m\^eme.
Cela montre la conjecture A. 

La conjecture $\Brm^-$ d\'ecoule de consid\'erations similaires. 
\end{proof}

\bigskip

\subsection{Groupe di\'edral infini}
Le m\^eme raisonnement que dans le cas fini permet de montrer facilement 
la proposition suivante~:

\bigskip

\begin{prop}
Si $S=\{s,t\}$ et $st$ est d'ordre infini dans $W$, alors les conjectures 
A et B sont vraies pour $(W,S)$. Les hyperplans essentiels de $(W,S)$ sont $\HC_s$, $\HC_t$, $\HC_{s-t}$ et 
$\HC_{s+t}$. 
\end{prop}

\bigskip

\subsection{Type ${\boldsymbol{F_4}}$}
Supposons dans cette sous-section, et dans cette sous-section seulement, 
que $(W,S)$ est de type $F_4$. Posons $S=\{s_1,s_2,t_1,t_2\}$ de 
telle sorte que le graphe de Coxeter de $(W,S)$ soit
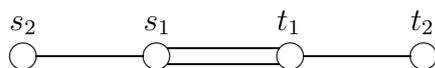
\begin{center}
\begin{picture}(160,28)
\put(0,10){\circle{10}}\put(-5,20){${s_2}$}
\put(50,10){\circle{10}}\put(45,20){${s_1}$}
\put(100,10){\circle{10}}\put(95,20){${t_1}$}
\put(150,10){\circle{10}}\put(145,20){${t_2}$}
\put(5,10){\line(1,0){40}}
\put(54,13){\line(1,0){42}}
\put(54,7){\line(1,0){42}}
\put(105,10){\line(1,0){40}}
\end{picture}
\end{center}
Notons $s=\{s_1,s_2\}$ et $t=\{t_1,t_2\}$, de sorte que $\Sba=\{s,t\}$. 

\bigskip

\begin{theo}[Geck]
Si $(W,S)$ est de type $F_4$, alors les conjectures A et $B^-$ sont vraies 
pour $(W,S)$. Avec les notations pr\'ec\'edentes, les hyperplans essentiels de 
$(W,S)$ sont $\HC_s$, $\HC_t$, $\HC_{s-2t}$, $\HC_{s-t}$, $\HC_{2s-t}$, 
$\HC_{s+2t}$, $\HC_{s+t}$ et $\HC_{2s+t}$.
\end{theo}

\bigskip

Le dessin de l'arrangement d'hyperplans du th\'eor\`eme pr\'ec\'edent 
(du moins, de son image inverse sous $\Pos$ dans $\RM[\Sba]^*$) est le suivant~:
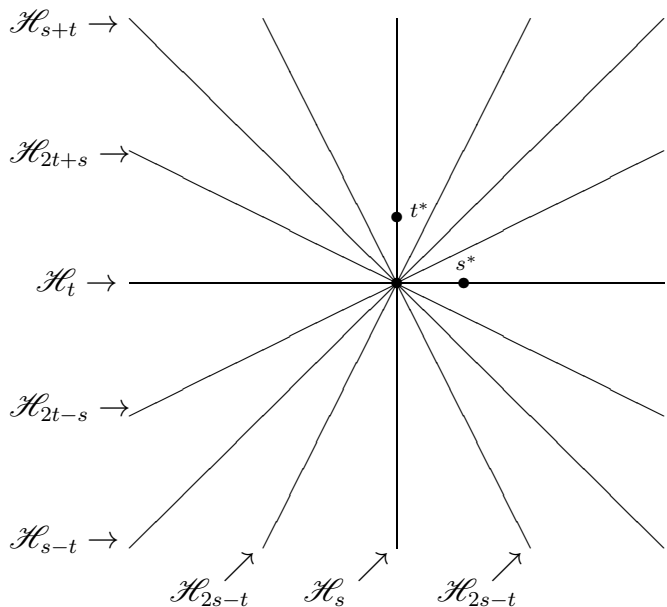
\bigskip

\begin{center}
\begin{picture}(200,200)
\put(0,100){\line(1,0){200}}
\put(0,0){\line(1,1){200}}
\put(100,0){\line(0,1){200}}
\put(0,200){\line(1,-1){200}}
\put(50,0){\line(1,2){100}}
\put(50,200){\line(1,-2){100}}
\put(0,50){\line(2,1){200}}
\put(0,150){\line(2,-1){200}}
\put(100,100){\circle*{4}}
\put(125,100){\circle*{4}}\put(122,105){$\SS{s^*}$}
\put(100,125){\circle*{4}}\put(105,125){$\SS{t^*}$}
\put(-35,97){$\HC_t \to$}
\put(-45,0){$\HC_{s-t} \to$}
\put(-45,195){$\HC_{s+t} \to$}
\put(-45,50){$\HC_{2t-s} \to$}
\put(-45,146){$\HC_{2t+s} \to$}
\put(65,-20){$\HC_s$}\put(85,-10){$\nearrow$}
\put(15,-20){$\HC_{2s-t}$}\put(35,-10){$\nearrow$}
\put(115,-20){$\HC_{2s-t}$}\put(135,-10){$\nearrow$}
\end{picture}
\end{center}

\vskip1cm

\noindent{\sc Remarque - } Tels quels, les r\'esultats de Geck \cite{geck f4} 
ne d\'emontrent pas compl\`etement les conjectures A et $\Brm^-$ en type $F_4$. 
Ses r\'esultats ne sont pas non plus \'enonc\'es sous cette forme. 
N\'eanmoins, comme le lecteur pourra le v\'erifier dans le sch\'ema 
de la preuve que nous donnons ci-dessous, Geck avait fait l'essentiel du travail: il a d\'etermin\'e les hyperplans essentiels. 
Nous avons de plus utilis\'e ses programmes {\tt GAP/CHEVIE} pour terminer 
le calcul lorsque certains param\`etres sont nuls.\finl

\bigskip

\begin{proof}
Dans \cite{geck f4}, M. Geck calcule les cellules de $W$ pour tous 
les choix de fonction $\ph$ telle que $\ph(S) \subseteq \G_{>0}$. 
Pour obtenir les cellules lorsque $\ph(S) \subseteq \G_{\geqslant 0}$, 
il suffit alors d'utiliser le corollaire \ref{KL semi-direct} 
(et la d\'ecomposition $W=\SG_3 \ltimes W(D_4)$, voir \ref{dyer bis}) 
ainsi que la connaissance des cellules pour un groupe de Weyl de 
type $D_4$, 
ce qui se calcule gr\^ace \`a {\tt GAP}~: nous l'avons fait en utilisant 
les programmes de M. Geck. Le corollaire \ref{symetrie} permet alors 
de d\'eterminer les cellules pour tous les choix de fonction $\ph$. 
Une discussion similaire \`a celle du cas di\'edral, tenant compte 
de tous les calculs effectu\'es, permet de v\'erifier les conjectures A et $\Brm^-$~: 
elle est simplifi\'ee par la pr\'esence de l'automorphisme de $W$ 
\'echangeant $s_i$ et $t_i$. 
\end{proof}

\subsection{Type ${\boldsymbol{\widetilde{G}_2}}$}
Supposons ici que $S=\{t,s_1,s_2\}$, que $(W,S)$ est un syst\`eme 
de Coxeter de type $\widetilde{G}_2$, et que le graphe de 
Coxeter est donn\'e par
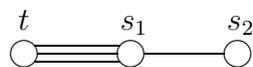
\begin{center}
\begin{picture}(100,30)
\put(10,10){\circle{10}}\put(8,19){$t$}
\put(50,10){\circle{10}}\put(46,19){$s_1$}
\put(90,10){\circle{10}}\put(86,19){$s_2$}
\put(55,10){\line(1,0){30}}
\put(15,10){\line(1,0){30}}
\put(14,13){\line(1,0){32}}
\put(14,7){\line(1,0){32}}
\end{picture}
\end{center}

\noindent Posons $s=\{s_1,s_2\}$ et identifions $t$ et $\{t\}$, de sorte 
que $\Sba=\{s,t\}$. Les diff\'erents calculs effectu\'es par J. Guilhot 
(voir \cite{guilhot} et \cite{guilhot2})
l'ont conduit a propos\'e la conjecture suivante~:

\bigskip

\bigskip

\noindent{\bf Conjecture G (Guilhot).} 
{\it Si $(W,S)$ est un groupe de Weyl affine de type $\widetilde{G}_2$, alors la 
conjecture A est vraie. Avec les notations pr\'ec\'edentes, les hyperplans 
essentiels sont $\HC_s$, $\HC_t$, $\HC_{s-t}$, $\HC_{s+t}$, 
$\HC_{2t-3s}$, $\HC_{2t+3s}$, $\HC_{t-2s}$ et $\HC_{t+2s}$.}

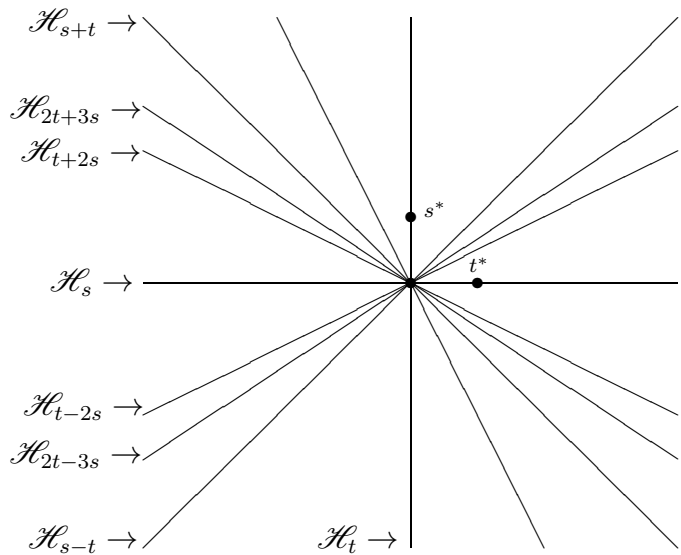
\bigskip

\begin{center}
\begin{picture}(250,200)
\put(50,100){\line(1,0){200}}
\put(50,0){\line(1,1){200}}
\put(150,0){\line(0,1){200}}
\put(50,200){\line(1,-1){200}}
\put(50,33.3){\line(3,2){200}}
\put(50,166.7){\line(3,-2){200}}
\put(100,200){\line(1,-2){100}}
\put(50,50){\line(2,1){200}}
\put(50,150){\line(2,-1){200}}
\put(150,100){\circle*{4}}
\put(175,100){\circle*{4}}\put(172,104){$\SS{t^*}$}
\put(150,125){\circle*{4}}\put(155,125){$\SS{s^*}$}
\put(15,97){$\HC_s \to$}
\put(7,0){$\HC_{s-t} \to$}
\put(7,195){$\HC_{s+t} \to$}
\put(5,50){$\HC_{t-2s} \to$}
\put(4,146){$\HC_{t+2s} \to$}
\put(115,0){$\HC_t \to$}
\put(0,162){$\HC_{2t+3s} \to$}
\put(0,32){$\HC_{2t-3s} \to$}
\end{picture}
\end{center}

\bigskip

\noindent{\sc Commentaire - } Le lecteur pourra se reporter \`a l'article 
\cite{bonnafe guilhot} pour des r\'esultats partiels r\'ecents corroborant 
la conjecture G.\finl

\bigskip

\subsection{Type ${\boldsymbol{B}}$\label{section B}}
Nous supposons ici que $(W,S)=(W_n,S_n)$, o\`u $W_n$ est un groupe 
de Weyl de type $B_n$ ($n \ge 2$) et $S_n=\{t,s_1,s_2,\dots,s_{n-1}\}$, 
de sorte que le graphe de Coxeter de $(W_n,S_n)$ soit
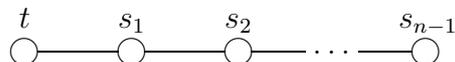
\begin{center}
\begin{picture}(170,40)
\put(10,10){\circle{10}}\put(8,19){$t$}
\put(50,10){\circle{10}}\put(45,19){$s_1$}
\put(90,10){\circle{10}}\put(85,19){$s_2$}
\put(160,10){\circle{10}}\put(150,19){$s_{n-1}$}
\put(15,10){\line(1,0){30}}
\put(55,10){\line(1,0){30}}
\put(95,10){\line(1,0){20}}
\put(155,10){\line(-1,0){20}}
\put(118,9.4){$\dots$}
\end{picture}
\end{center}
Identifions $t$ et $\{t\}$ et posons $s=\{s_1,s_2,\dots,s_{n-1}\}$. 
On a ainsi $\Sba=\{s,t\}$. Dans ce cas, les conjectures de Geck, Iancu, 
Lam et l'auteur \cite[Conjectures A et B]{bgil} sugg\`erent la 
conjecture suivante~:

\bigskip

\noindent{\bf Conjecture C.} 
{\it Les conjectures A et B sont vraies pour $(W_n,S_n)$. Les hyperplans 
essentiels sont $\HC_s$, $\HC_t$, $\HC_{t-is}$ ($1 \le i \le n-1$) et 
$\HC_{t+is}$ ($1 \le i \le n-1$).}

\bigskip

\noindent{\sc Commentaire - } Le lecteur pourra se reporter aux articles  
\cite{lacriced}, \cite{bonnafe two}, \cite{bgil}, \cite{bonnafe B} et \cite{bonnafe guilhot} (dans l'ordre chronologique) 
pour des r\'esultats corroborant la conjecture C.\finl


%
%
%

\end{document}